\documentstyle[12pt]{article}     
\def\m@th{\mathsurround=0pt }
\def\eqalign#1{\null\,\vcenter{\openup\jot \m@th
   \ialign{\strut\hfil$\displaystyle{##}$&$
      \displaystyle{{}##}$\hfil \crcr#1\crcr}}\,}
\newtheorem{theorem}{Theorem}[section]

\newtheorem{corollary}[theorem]{Corollary}
\newtheorem{lemma}[theorem]{Lemma}

\def\ad{{\rm ad }\ }
\def\adr{{\rm ad}_r\ }

\def\wt{{\rm wt }\ }

\font\goth=eufm10

\def\gg{\hbox{\goth g}}

\def\gn{\hbox{\goth n}}

\def\gh{\hbox{\goth h}}

\begin{document}

\begin{center}
{\large{\bf  Quantum Symmetric Pairs and Their Zonal Spherical 
Functions } }
\vskip 5mm Gail Letzter\footnote{supported by NSA grant no. 
MDA904-01-1-0033.

AMS subject classification 17B37}

Mathematics Department

Virginia Polytechnic Institute

and State University

Blacksburg, VA 24061

\end{center}
\medskip
\begin{abstract}
 We study   the space of biinvariants and   zonal 
spherical functions associated to quantum symmetric pairs in the 
maximally split case.
Under the obvious restriction map,    the space of biinvariants  
is proved isomorphic to the   Weyl group invariants of  the   character group ring 
 associated to the restricted roots.
  As a consequence,  
 there is either a unique set, or an (almost) unique  two-parameter 
 set of   Weyl group invariant quantum zonal spherical 
functions
associated to an irreducible symmetric pair. Included is  a 
complete and explicit list of the generators and relations for the left coideal 
subalgebras of the quantized enveloping algebra  used to form 
quantum  symmetric   pairs. 
 
\end{abstract}

\centerline{INTRODUCTION}
\medskip

   The representation theory of semisimple Lie algebras and Lie groups 
has been closely intertwined with the theory of symmetric spaces 
since E. Cartan's pioneering work in the 1920's.     
A beautiful classical result shows that 
 the zonal spherical functions of these spaces can be identified with a family of 
 orthogonal polynomials.
 With the 
introduction of quantum groups in the 1980's, it was natural to look 
for and study quantum versions of symmetric spaces. This search 
became especially compelling as $q$ orthogonal polynomials, which are 
obvious candidates for 
quantum zonal spherical functions, appeared in the literature
(see for example [Ma] and [K2]). 
However,  the theory of quantum symmetric spaces was 
initially slow to develop because it was not obvious how to form 
them.

Let $\gg$ be a complex semisimple Lie algebra and let $\theta$ be an 
involution of $\gg$.   A classical (infinitesimal) symmetric pair 
consists of the Lie algebra $\gg$ and the fixed Lie subalgebra 
$\gg^{\theta}$.  
In his fundamental paper [K1], Koornwinder was the first to use an 
infinitesimal approach to  construct 
quantum symmetric spaces.  
In particular, he found
quantum analogs of a symmetric pair of Lie 
algebras for $\gg={\bf sl}\ 2$ using twisted primitive elements.
As this method was generalized to other examples, it became clear 
that 
coideal analogs of $U(\gg^{\theta})$ inside the 
quantized enveloping algebra $U_q(\gg)$ had the potential to produce a 
``good'' theory of quantum symmetric spaces (see for example [N],[NS], and [DN].) In his 
1996
survey paper [Di], Dijkhuizen describes  the philosophy of this 
approach
in three major steps which we briefly rephrase here.

\begin{enumerate}
\item[(1)] Find coideals inside the quantized enveloping algebra 
which can be used to form quantum symmetric pairs.  Show that the  
finite-dimensional spherical modules can be parametrized using the 
dominant integral restricted roots.
\item[(2)] Determine the image of the 
space of biinvariant functions associated to a quantum symmetric pair
 inside the character ring 
of the root system of $\gg$.
\item[(3)] Realize zonal spherical functions as $q$-orthogonal 
polynomials by computing the radial components of appropriate central 
elements of the quantized enveloping algebra.
\end{enumerate}

   This  paper is one of a series  by the author 
   intended to
answer these questions when the triangular decomposition of $\gg$ has 
been chosen so that the Cartan subalgebra is maximally split with 
respect to $\theta$.     Quantum analogs  of the enveloping algebra
$U(\gg^{\theta})$   are constructed in [L2] and [L4] by describing 
their   generators. These subalgebras of $U_q(\gg)$ are 
further characterized as the  unique maximal left coideal subalgebras   which specialize 
to $U(\gg^{\theta})$ as $q$ goes to $1$.  The finite-dimensional 
spherical modules associated to quantum symmetric pairs are classified in [L3],
thus completing problem (1).

The   main result of this paper is  a comprehensive  answer to problem 
(2).
Let ${\cal B}_{\theta}$ denote the set  of  
coideal subalgebras associated to  $\gg,\gg^{\theta}$.
For each pair of subalgebras  $B$ and $B'$ in ${\cal B}_{\theta}$, 
the space of biinvariants  ${ }_{B'}{\cal H}_{B}$  
consists of the  left $B'$ and right $B$ invariants 
inside the quantized function algebra $R_q[G]$ 
corresponding to $U_q(\gg)$. The vector space ${ }_{B'}{\cal H}_{B}$   
is an algebra and can be written as a direct sum of eigenspaces with 
respect to the action of the center of $U_q(\gg)$.  The eigenvectors 
with respect to this action are called zonal spherical functions.
 
Elements of the quantized function algebra can be thought of 
as functions on $U_q(\gg)$. Recall that $U_q(\gg)$ contains a 
multiplicative subgroup $T$ corresponding to the root lattice of the 
root system of $\gg$.  Let $\Sigma$ denote the restricted root system 
associated to $\gg,\theta$ and let ${\cal C}[2\Sigma]$ be the 
character group ring of the weight lattice corresponding to $2\Sigma$. 
   We prove that restriction to $T$ of the space of biinvariants 
   (considered as functions on $U_q(\gg)$) induces  an 
injection   $\Upsilon$  from   ${ }_{B'}{\cal H}_{B}$
into ${\cal C}[2\Sigma]$.

Let $W$ denote the restricted Weyl group associated to 
$\Sigma$. Let ${\bf H}$ denote the group of Hopf algebra automorphisms of $U_q(\gg)$ 
which fix elements of $T$. Note that  ${\bf H}$  acts on ${\cal 
B}_{\theta}$ by sending an algebra to its image 
 under the automorphism. Given $B\in {\cal B}_{\theta}$
 and any ${\bf H}$ orbit ${\cal O}$ of ${\cal 
B}_{\theta}$  we find 
$B'\in {\cal O}$ such that  $\Upsilon(_{B'}{\cal H}_{B})$ is $W$ 
invariant.  We then show, in this case, that $\Upsilon(_{B'}{\cal H}_{B})$ 
  is
equal to ${\cal C}[2\Sigma]^W$.
 
The action of ${\bf H}$ on ${\cal B}_{\theta}$ can be extended 
componentwise to an 
action of  the group ${\bf H}\times {\bf H}$  on ${\cal 
B}_{\theta}\times{\cal B}_{\theta}$.  We show that the images 
 under $\Upsilon$ of two different spaces 
 of biinvariants associated to  pairs in the same ${\bf H}\times {\bf H}$
 orbit  are related by an automorphism 
of ${\cal C}[2\Sigma]$.  Thus the image of any space of biinvariants 
under $\Upsilon$
is  just a
translation of
 ${\cal C}[2\Sigma]^W$ via an automorphism of ${\cal C}[2\Sigma]$.
   Furthermore, this automorphism  restricts 
to an eigenvalue preserving function on the corresponding sets of  zonal 
spherical functions.

 A zonal spherical family associated to 
a ${\bf H}\times {\bf H}$ orbit 
 is  the image under $\Upsilon$  of a specially chosen basis of zonal spherical 
 functions in  ${  }_{B'}{\cal 
H}_{B}$  for some pair 
$(B,B')$ in this orbit. We
show that there is a vitually unique    $W$ invariant zonal spherical family 
associated to each  ${\bf H}\times {\bf H}$ orbit  of ${\cal 
B}_{\theta}\times {\cal B}_{\theta}$.  The obstruction to uniqueness,
a small finite subgroup of ${\rm Aut}{\cal C}[2\Sigma]^W$, is 
just the trivial group for 
most choices  of $\gg$ and $\theta$.

The pair $\gg,
\gg^{\theta}$,  or more precisely, $\gg,\theta$, is called irreducible if $\gg$ cannot be written as 
the direct sum of two semisimple Lie subalgebras which both admit $\theta$ 
as an involution. In [A] (see also [He, Chapter X, Section F]), Araki gave a complete list of all the irreducible pairs 
$\gg,\theta$. In the final section of the paper,
using this list and the construction ([L4, Section 7]) of the algebras 
in ${\cal B}_{\theta}$, we   explicitly describe the 
generators and relations for every coideal subalgebra in ${\cal 
B}_{\theta}$ associated to 
each irreducible pair $\gg,{\theta}$. It follows from this 
classification that if  $\gg,\theta$ is 
irreducible, either ${\cal B}_{\theta}\times {\cal 
B}_{\theta}$ consists of a single ${\bf H}\times {\bf H}$ orbit or a 
two-parameter set of orbits.  In the first case, there is a 
unique $W$ invariant zonal spherical family associated to ${\cal 
B}_{\theta}$.  For the latter case, there is either a unique 
two-parameter set of 
$W$ invariant zonal spherical families associated to ${\cal 
B}_{\theta}$ or exactly 
two two-parameter sets.

In a future paper, we will address problem (3) and 
realize  these $W$ invariant zonal spherical 
families as   sets of $q$ hypergeometric polynomials. Apparently,
when $\Sigma$ is reduced and ${\cal B}_{\theta}$ is a single ${\bf 
H}\times {\bf H}$ orbit,  the zonal spherical functions are Macdonald polynomials. 
Here, one of the parameters is a power of the other.  This would generalize the results in [N] as well as 
those announced in [NS] 
concerning  quantum symmetric pairs which can be constructed 
using solutions to reflection equations. Other $q$ 
hypergeometric functions such as Askey-Wilson 
polynomials with two indeterminates arise in the remaining cases.  This 
should remind the reader of   the zonal spherical functions computed on the 
quantum $2$-sphere  in [K1] and 
 on the family of quantum projective spaces 
studied in [DN].

The remainder of the paper is organized as follows. Section 1 sets notation and 
reviews basic facts about involutions, restricted root systems,
  and quantized enveloping algebras.  In 
Section 2, we review the construction of the quantum analogs of 
$U(\gg^{\theta})$ inside $U_q(\gg)$.   Section 3 is a 
study of spherical modules and their spherical vectors.   
Fine information about the possible weights of weight vectors which show up as 
summands of spherical vectors is obtained.   The next three sections of 
the paper are
devoted to the space of biinvariants  and their 
zonal spherical functions.  In 
Section 4, the map $\Upsilon$ restricted to a space of biinvariants
 is shown to be injective with image contained in ${\cal C}[2\Sigma]$.  Section 5 
establishes a criterion for choosing $B$ and $B'$ so that ${}_{B'}{\cal 
H}_{B}$ is $W$ invariant.  Section 6 is a study of zonal spherical 
families.   We show how to relate them using elements of the group 
${\bf H}\times {\bf H}$ and determine necessary and sufficient 
conditions for a zonal spherical family to be $W$ invariant. In 
Section 7, we give a complete and explicit list of 
generators and relations for the algebras in ${\cal B}_{\theta}$ 
associated to 
each irreducible pair $\gg, {\theta}$.

\section{Background and Notation}
Let ${\bf C}$ denote the complex numbers, ${\bf R}$ denote the real 
numbers, ${\bf Q}$ denote the rational numbers, ${\bf Z}$ denote the 
integers, and ${\bf N}$ denote the nonnegative integers.  
Let ${\gg}={\gn}^-\oplus {\gh}\oplus {\gn}^+$ be a 
semisimple Lie algebra of rank $n$ over  ${\bf C}$ with Cartan matrix $(a_{ij})$. Write $\Delta$  
for the set of 
roots  of ${\gg}$ and let $\pi=\{\alpha_1,\dots, 
\alpha_n\}$ be a fixed set of positive simple roots.
 Set $Q(\pi)=\sum_{1\leq 
i\leq n}{\bf Z}\alpha_i$ and $Q^+(\pi)=\sum_{1\leq i\leq n}{\bf N}\alpha_i$. 
 Let $P^+(\pi)$ 
denote the set of dominant integral weights associated to the root 
system $\Delta$ and  $P(\pi)$ denote the 
weight lattice.  Write $(\ , \ )$ for 
the Cartan inner product. Let $\leq$ denote the standard partial order 
on $Q(\pi)$.   In particular, given two elements $\lambda$ and $\beta$ 
in $Q(\pi)$, we say that $\lambda\leq\beta$ if and only if $\lambda-\beta\in 
Q^+(\pi)$.

Let $\theta$ be a maximally split  Lie algebra  involution of ${\bf g}$ 
with respect to the fixed Cartan subalgebra ${\gh}$ and triangular 
decomposition of $\gg$ 
in the sense of [L4, Section 7]. 
Write ${\gg}^{\theta}$ for the corresponding fixed Lie subalgebra. 
The 
involution $\theta$ induces an involution $\Theta$ of the root system 
of ${\gg}$ and thus an automorphism of $\gh^*$.  Let 
$\pi_{\Theta}$ be the set of simple roots 
$\{\alpha_i|\Theta(\alpha_i)=\alpha_i\}$ 
and write $\Delta_{\Theta}$ for the corresponding root system generated by 
$\pi_{\Theta}$. Let $p$ be the permutation on $\{1,\dots,n\}$ 
corresponding to a diagram automorphism of $\pi$  such that
$\Theta(\alpha_i)+\alpha_{p(i)}\in \sum_{\alpha_i\in \pi_{\Theta}}{\bf 
Z}\alpha_i$ for each 
$\alpha_i\in \pi-\pi_{\Theta}$. (For more information about 
involutions, see for example [D, 1.13] and  [L4, Section 7]).

Given $\beta\in \gh^*$, set $\tilde\beta=  1/2(\beta-\Theta(\beta))$.
 Let $\Sigma$ denote the restricted root system associated to 
$\gg,\theta$ ([Kn, Chapter VI,  Section 4], [He, Chapter X, Section 
F], or [L4, Section 7]). 
Note that 
$\Sigma$ can be identified with the set
$\{\tilde\alpha|\alpha\in\Delta-\Delta_{\Theta}\}$
using the Cartan inner product as the inner product for  $\Sigma$.
Moreover,
$\{\tilde\alpha_i| 
\alpha_i\in\pi-\pi_{\Theta}\}$
is a set of positive simple roots for the root system $\Sigma$. 
Set $Q(\Sigma)$ equal to the root lattice 
$\sum_i{\bf Z}\tilde\alpha_i$ 
of $\Sigma$ and
set $P(\Sigma)$  equal to the weight lattice associated to the root 
system $\Sigma $. Let $P^+(\Sigma)$ denote the subset of 
$P(\Sigma)$ 
consisting of dominant integral weights.  
Set $P(2\Sigma)=\{2 \lambda| \lambda\in P(\Sigma)\}$.

Let $q$ be an indeterminate and set $q_i=q^{(\alpha_i,\alpha_i)/2}$ 
for each $1\leq i\leq n$.
Write $U= U_q({\gg})$  for the quantized enveloping algebra generated
by $x_i$, $y_i$, $t_i^{\pm 1}$, $1\leq i\leq n$, over 
the algebraic closure ${\cal C}$ of   ${\bf C}(q)$ (See  [L4, 
Section 1, (1.4)-(1.10)] or [Jo, 3.2.9] for 
 relations.) Let $U^-$ be the subalgebra of $U$ generated by $y_i,1\leq i\leq n$ 
 and let $U^+$ be the subalgebra of $U$ generated by 
 $x_i,1\leq i\leq n$.  Given an  
integral domain $D$ containing ${\cal C}$, we  set 
$U_{D}=U\otimes_{{\cal C}}{D}$.
 
 The algebra $U$ is a Hopf algebra with 
 comultiplication map
 ${\it\Delta}$, antipode $\sigma$, and counit $\epsilon$.
 Let $U_+$
denote the augmentation ideal of
$U$, which is the kernel of the counit  map $\epsilon.$   Given a subalgebra $A$ of $U$, we write
$A_+$ for the intersection of 
$A$ with $U_+$.

 We use 
 Sweedler notation for the coproduct.   In particular, we write
 $${\it\Delta}(a)=\sum a_{(1)}\otimes a_{(2)}$$ for each $a\in  U$. 
Recall that  a subalgebra $S$ of $U$ is called a left coideal 
subalgebra if 
 ${\it\Delta}(a)\in U\otimes S$ for all $a\in S$.  
 
 Let
 ${\rm ad}_r$ denote the right adjoint action defined by $$(\adr 
 a)b=\sum\sigma(a_{(1)})ba_{(2)}\leqno{(1.1)}$$ for all $a\in U$ and $b\in U$.  For 
 example, $$(\adr 
 x_j)b=-t_j^{-1}x_jb+t_j^{-1}bx_j{\rm \quad and \quad }(\adr 
 y_j)b=by_j-y_jt_jbt_j^{-1}\leqno(1.2) $$ for  $b\in U$ and $1\leq j\leq n$.

Let $T$ be the group generated 
by the $t_i$, $1\leq i\leq n$ and let $\tau$ be the 
isomorphism from $Q(\pi)$ to $T$ which sends $\alpha_i$ to $t_i$. 
Write    $U^o$ for the group algebra generated by $T$.
Define the subgroup $T_{\Theta}$  of $T$ by
 $$T_{\Theta}=\{\tau(\lambda)|\lambda\in Q(\pi) {\rm \ and \ }
\Theta(\lambda)=\lambda\}.$$

Let $M$ be a $U$ module. We say that a nonzero vector $v\in U$ has 
weight $\lambda\in \gh^*$ provided that $\tau(\lambda)\cdot 
v=q^{(\gamma,\lambda)}v$ for all $\tau(\lambda)\in T$. 
For any vector subspace $V\subset M$, write $V_{\gamma}$ for the subspace of $V$ spanned by the $\gamma$
weight vectors.  If $V$ is a vector subspace of $U$, then $V_{\gamma}$ is the subspace 
of $V$ consisting of $\gamma$ weight vectors with respect to the 
adjoint action. Now suppose that  
$F$ is both a subalgebra   and $\ad T$ submodule of $U$.  Then for any 
subgroup $S$ of $T$, we write  $FS$ for the subalgebra of $U$ 
generated by $F$ and $S$.  Note that $FS$ is 
 spanned as a vector space over ${\cal C}$ by elements $as$ with $a\in F$ and $s\in S$.

\medskip
\section{Quantum Symmetric Pairs}
Quantum analogs of the pair $U(\gg),U(\gg^{\theta})$ associated to a 
maximally split involution 
 are 
constructed and characterized up to Hopf algebra automorphism 
in [L2] and [L4, Section 7]. These analogs consist of  $U$ and   a maximal left 
coideal subalgebra $B$ of $U$ which specializes to $U(\gg^{\theta})$ as 
$q$ goes to $1$  (The reader is referred to [L4, Section 1, 
following (1.10)] for the precise definition of specialization used
here, [L4, (7.25)] for the notion  of  maximal,  and   [L4, Theorem 7.5].) 
 We review the construction of these subalgebras here.

From now on, $\Theta$ will denote an involution of the root system $\Delta$ induced by a 
maximally split involution of $\gg$.     Let  
${\cal M} $ be the subalgebra of $U$ generated by  $x_i$, 
$y_i$, $t_i^{\pm 1}$
for $\alpha_i\in \pi_{\Theta}$.   Set ${\cal 
M}^+$ 
equal to the subalgebra of ${\cal M}$ generated by the 
$x_i$, $\alpha_i\in\pi_{\Theta}$.  

Let $[\Theta]$ be the 
set of all maximally split involutions of $\gg$ which induce the 
involution $\Theta$ on $\Delta$. Note that two involutions in 
$[\Theta]$ are conjugate to each other via a Lie algebra 
automorphism of $\gg$ which fixes the Cartan subalgebra. 
Let $\theta$ be the particular choice of  maximally split involution 
 in $[\Theta]$ 
 as described in [L4, Section 7, the discussion following (7.5) ] 
 and $\tilde\theta$ be the  automorphism of $U$ 
defined in [L4, Theorem 7.1] which specializes to $\theta$.
We describe the action of $\tilde\theta$ on $y_i$ for 
$\alpha_i\notin\pi_{\Theta}$. (Note that replacing each $x_j$ by $e_j$, 
each $t_j$ with $1$, and setting $q=1$ will recover the action of 
$\theta$ on $f_i$.  Since $\theta$ 
is the identity on the positive and negative root vectors 
corresponding to roots in  $\pi_{\Theta}$ and the action of $\theta$ on 
$\gh$ corresponds to the action of $\Theta$ on $\gh^*$,  this information is 
enough to determine $\theta$.)
 
Let $\pi^*$ be the  subset of $\pi-\pi_{\Theta}$ consisting of 
all $\alpha_i\notin\pi_{\Theta}$ such that either $i=p(i)$ or 
$i<p(i)$.  Write $x_i^{(m)}$ and $y_i^{(m)}$ where $m\in {\bf N}$
 for the divided powers  of $x_i$ and 
$y_i$ repectively  (see 
[Jo, 4.3.14 and 1.2.12]).
Given $i$ such that $\alpha_i\in\pi^*$, 
there exists a sequence $\alpha_{i_1},\dots, \alpha_{i_r}$ consisting of 
elements in $\pi_{\Theta}$ and positive integers $m_1,\dots, m_r$ 
subject to  
$$\tilde\theta(y_i)=(\adr x_{i_1}^{(m_1)})\cdots 
(\adr x_{i_r}^{(m_r)})t_{p(i)}^{-1}x_{p(i)}\leqno(2.1)$$
and $$\tilde\theta(y_{p(i)})=(-1)^{(m_1+\cdots+m_r)}(\adr 
x_{i_r}^{(m_r)})\cdots 
(\adr x_{i_1}^{(m_1)})t_{i}^{-1}x_{i}.\leqno(2.2)$$  Using (1.1), it is 
straightforward to check that $(\adr y_j)t_i^{-1}x_i=0$ whenever 
$i\neq j$.   Thus $t_i^{-1}x_i$ is a lowest weight vector with respect 
to the action of $\adr{\cal M}$
for all $\alpha_i\notin\pi_{\Theta}$. It follows  as in the 
classical case ([L4, discussion following (7.6)]) that these sequences 
satisfy the following property.  For each $r\geq s> 0$, 
$$X_s=(\adr x_{i_{s}}^{(m_{s})})(\adr x_{i_{s+1}}^{(m_{s+1})})\cdots 
(\adr x_{i_r}^{(m_r)})t_{p(i)}^{-1}x_{p(i)}\leqno(2.3)$$ is 
a highest weight vector for the action of $\adr x_{i_s}$ and $X_{s+1}$ is a 
lowest weight vector   for the action 
of $\adr y_{i_s}$.    
  Moreover, $\tilde\theta(y_i)$ is a highest 
weight vector for the action of $\adr{{\cal M}}$.

  Set ${\cal S}$ equal to the 
 subset  of $\pi^*$ consisting of  $\alpha_i$ such that 
$\Theta(\alpha_i)=-\alpha_i$ and $2(\alpha_i,\alpha_j)/(\alpha_j,\alpha_j)$ 
is even
for all $\alpha_j$ such that $\Theta(\alpha_j)=-\alpha_j$. Let ${\cal D}$ 
denote the subset of $\pi^*$ consisting of those $\alpha_i$ 
such that $i\neq p(i)$ and $(\alpha_i,\Theta(\alpha_i))\neq 0$. 
Let ${\cal Z}$ be an integral domain containing  the complex numbers 
${\bf C}$ and 
contained in some field extension ${\cal C}({\cal Z})$ of ${\cal C}$.
Write ${\cal Z}^{\times}$ for the nonzero elements of ${\cal Z}$. 
 Define 
$${\bf S}({\cal Z})=\{{\bf s}\in {\cal Z}^{n}|s_i=0{\rm\ for\ } 
\alpha_i\notin{\cal S}\}$$
and 
$${\bf D}({\cal Z})=\{{\bf d}\in ({\cal Z}^{\times})^{n}|d_i=1{\rm\ for\ } 
\alpha_i\notin{\cal S}\}.$$ 

   Given ${\bf s}\in {\bf S}({\cal Z})$ and 
${\bf d}\in {\bf S}({\cal Z})$, the subalgebra
$ B_{\theta,{\bf s},{\bf d}}$ of $U_{{\cal C}({\cal Z})}$ is  generated by 
$T_{\Theta}$, ${\cal M}$, and the elements 
$B_{i,s_i,d_i}^{\theta}$ for $\alpha_i\in \pi-\pi_{\Theta}$ where  
$$B_{i,s_i,d_i}^{\theta}= 
y_it_i+d_i\tilde\theta(y_i)t_i+s_it_i.\leqno{(2.4)}$$
 We     abbreviate
 $B_{i,s_i,d_i}^{\theta}$ by $B_i$ when $s_i, d_i, $ and 
 $\tilde\theta$ can be understood from the context.  Note that if $i=p(i)$, then 
the definitions of  $B_i$ and $B_{p(i)}$ using (2.4) agree.

  The next theorem follows as in the  proof of [L4, Theorem 7.2]. (See 
  also [L4, Variations 1 and 2].)

 \begin{theorem} Let ${\cal Z}$ be an extension field ${\cal C}$.
  For each ${\bf s}\in {\bf S}({\cal Z})$ and ${\bf d}\in {\bf 
 D}({\cal Z})$, the algebra $B_{{\theta},{\bf s},{\bf d}}$ is a 
 left coideal subalgebra of $U_{{\cal Z}}$.
\end{theorem}

\medskip

Recall that a pair $\gg,\theta$ is called ireducible if $\gg$ cannot 
be written as the direct sum of two semisimple Lie subalgebras which 
both admit $\theta$ as an involution. 
In Section 7, we give
a   complete list  of   coideal subalgebras 
$B_{\theta,{\bf s},{\bf d}}$ associated to  
  irreducible pairs 
$\gg,\theta$. 
From this list, it can be seen that  if $\gg,\theta$ is irreducible, 
then either both ${\cal S}$ and ${\cal D}$ 
are the empty set, or one is the empty set and the other consists of 
exactly one root.

  Set ${\cal B}_{\theta}$ equal to the union of the orbits of all
    $B_{{\theta},{\bf 
s},{\bf d}}$ under the action of ${\bf H}$ where
${\bf s}\in{\bf S}({\cal C})$ and ${\bf d}\in {\bf 
D}({\cal C})$. More generally, given an integral domain  ${\cal Z}$ 
containing
${\cal C}$, let ${\cal B}_{\theta}({\cal Z})$ be the set of 
algebras $B$ which are isomorphic via a Hopf 
algebra automorphism fixing  
${\cal M} T_{\Theta}$ of  $U_{{\cal Z}}$ to $B_{{\theta},{\bf 
s},{\bf d}}$, for some  ${\bf s}\in{\bf S}({\cal Z})$ and ${\bf d}\in {\bf 
D}({\cal Z})$.  
Note that if both ${\cal S}$ and ${\cal D}$ are empty, 
 then  ${\cal B}_{\theta}$   consists of exactly one 
 orbit under this action. It follows from the description of the group 
 of  Hopf 
 algebra automorphisms of $U$ in [Jo, Section 10.4] that the set $\{B_{\theta,{\bf 
 s},{\bf d}}|\ {\bf s}\in {\bf S}({\cal C}){\rm \ and \ }{\bf d}\in {\bf 
 D}({\cal C})\}$ is a complete set of distinct 
 representatives for the orbits of ${\cal B}_{\theta}$ under ${\bf H}$. 
 
 Let $\hat U$ denote the ${\bf C}[q]_{q-1}$ subalgebra of $U$ 
generated by $x_i,y_i,t_i^{\pm 
1},$ 
and $(t_i-1)/(q-1)$ for $1\leq i\leq n$.  Let ${\cal B}_{\theta,1}$ 
be the set of maximal left coideal subalgebras of $U$ such that 
  $B\cap\hat U$ specializes to $U(\gg^{\theta})$ 
as $q$ goes to $1$.  In particular,  ${\cal B}_{ \theta,1}$ is the set of quantum 
analogs of $U(\gg^{\theta})$ in the sense of [L2] and [L4].  By [L4, 
Theorem 7.5], every algebra in ${\cal B}_{\theta,1}$ is isomorphic to 
an appropriate $B_{{\theta},{\bf s},{\bf d}}$ for some ${\bf s}\in 
 {\bf S}({{\bf C}[q]}_{(q-1)})$ and ${\bf d}\in {\bf 
 D}({{\bf C}[q]}_{(q-1)})$ via some Hopf algebra automorphism in ${\bf H}$.

It should be noted that for certain $\theta$, 
some of the  algebras in ${\cal B}_{\theta}$ have appeared in the 
literature in papers not written by the author.
The most prominent of these algebras is the nonstandard
quantum analog $U'_q({\bf so}\ n)$ of $U({\bf so}\ n)$ introduced by  Gavrilik and 
 Klimyk
in [GK].  This algebra  is equal to 
$B_{\theta,{\bf s},{\bf 
d}}$ when $\gg,\theta$ is of type {\bf AI} as described in Section 7.    The representation 
theory of $U'_q({\bf so}\ n)$ has been extensively investigated in a series of papers by 
Gavrilik, Klimyk,  and others (see for example [GI], [GIK], 
and [HKP]).   In [NS], Noumi and Sugitani 
construct subalgebras $U^{tw}_q(\gg^{\theta})$ of $U$ which are 
quantum analogs  of $U(\gg^{\theta})$ 
 for many of 
the pairs 
$\gg,\theta$, when  $\gg$ is of classical type.  Their approach is 
completely different and involves solutions to reflection equations.  
Nevertheless, it is shown in [L2, Section 6], that for a given pair 
$\gg,\theta$, each subalgebra 
$U^{tw}_q(\gg^{\theta})$ constructed in [NS] belongs to ${\cal B}_{\theta}$.

 \medskip
 \section{Spherical modules}
 
Let $L(\lambda)$ 
denote the finite-dimensional simple $U$ module with highest weight 
$\lambda$ where $\lambda\in P^+(\pi)$. Given a left coideal 
subalgebra $B$  in ${\cal B}_{\theta}$ of $U$, we say that $L(\lambda)$ is 
spherical with respect to $B$  if  $L(\lambda)^B=\{v\in 
L(\lambda)|bv=\epsilon(b)v\}$ is a one-dimensional space. Necessary 
and sufficient conditions for $L(\lambda)$ to be spherical with 
respect to  
elements $B$ in ${\cal B}_{\theta,1}$ are established in [L3, 
Section 4]. This result is extended
in this paper to all algebras in ${\cal B}_{\theta}$.

Let $\{\zeta_i|\alpha_i\in {\cal S}\}\cup\{\eta_j|\alpha_j\in {\cal 
D}\}$ 
be a set of $|{\cal S}|+|{\cal D}|$ indeterminates.   Note that $-1$ cannot be written as a sum of squares 
in the field  ${\bf R}(q)(\zeta,\eta)$  generated by these indeterminates over ${\bf R}(q)$.
It follows that ${\bf R}(q)(\zeta,\eta)$ is a formally real field.  Thus 
we can choose a set of positive elements  for the field ${\bf 
R}(q)(\zeta,\eta)$  containing   $\eta_j$, for each 
$\alpha_j\in{\cal D}$. Let ${\cal K}$ denote the algebraic closure 
of ${\bf R}(q)(\zeta,\eta)$.
Set ${\cal R}$ equal to the real algebraic closure  of ${\bf 
R}(q)(\zeta,\eta)$ inside 
of ${\cal K}$ and let ${\cal R}^+$ denote the set of positive 
elements of ${\cal R}$. 
  (For more information on formally real 
fields and real algebraic closures, the reader is referred to [J, 
Chapter 11].)

Note that ${\cal K}={\cal R}+i{\cal R}$ and thus  complex conjugation extends 
to ${\cal K}$.
Let $\kappa$ be the conjugate linear antiautomorphism 
of   $U_{\cal K}$     defined by
 $\kappa(x_i)=y_it_i$,
$\kappa(y_i)=t_i^{-1}x_i$ and $\kappa(t)=t$ for all $t\in T$. 
In particular, $\kappa$ restricts to an antiautomorphism of $U_{\cal 
R}$ and acts as conjugation on elements of ${\cal C}$.  Note 
that $\kappa^2$ is just the identity. 
Let ${\bf H}_{\cal R}$ denote the subgroup of ${\bf H}$ whose elements restrict to Hopf algebra 
automorphisms of $U_{\cal R}$.

Let ${\bf s}(\zeta)=(s(\zeta_1),\dots, s(\zeta_n))$ be the element of ${\bf S}({\cal K})$ such that 
$s(\zeta_i)=\zeta_i$ for all $\alpha_i\in {\cal S}$.  Similarly, 
let $  {\bf d}(\eta)=(d(\eta_1),\dots, d(\eta_n))$ be the element 
of ${\bf D}({\cal K})$ such that $d(\eta_j)=\eta_j$ for all 
$\alpha_j\in {\cal D}$.  Since $d(\eta_j)=1$ for 
$\alpha_j\notin{\cal D}$, it follows  that $d(\eta_j)\in {\cal 
R}^+$ for all $1\leq j\leq n$.
By the definition of real algebraic closures,   
every positive element in ${\cal R}$ has a square root 
in ${\cal R}$. Let $ d(\eta_i)^{1/2}$ denote the positive square root of 
$  d(\eta_i)$ for 
each $i$.
We have the following 
generalizations of [L3, Lemma 3.2] and [L4, Theorem 7.6].

\begin{theorem}  There exists a Hopf algebra automorphism 
$\varphi\in {\bf H}_{\cal R}$ such that 
$\varphi\kappa\varphi^{-1}(B_{\theta,{\bf s}(\zeta),{\bf 
  d}(\eta)})=B_{\theta,{\bf   s}(\zeta),{\bf    d}(\eta)}$.
Moreover $B_{{\theta},{\bf 
  s}(\zeta),{\bf  d}(\eta)}$ acts semisimply on all finite-dimensional 
$U_{\cal K}$-modules.
 \end{theorem}
 
 \noindent
 {\bf Proof}: The first assertion implies the second statement  by 
 [L4, Section 2]. The 
 proof presented here for the first assertion 
 does not follow  the original argument given in [L3] which uses
 specialization at $q=1$.
 Instead, it is based closely on   the discussion leading up to
  [L4, Theorem 7.6].

 Set $B=B_{\theta,{\bf  s}(\zeta),{\bf 
  d}(\eta)}$ and $B_i=B^{\theta}_{i,s(\zeta_i), 
 d(\eta_i)}$ for each $\alpha_i\notin\pi_{\Theta}$.  Note  that $$(\varphi\kappa\varphi^{-1})({\cal 
M}T_{\Theta})={\cal M}T_{\Theta}$$ and $(\varphi\kappa 
\varphi^{-1})^2=\kappa^2$ is  the identity on $U_{\cal K}$ for all 
$\varphi\in {\bf H}_{\cal R}$.
Hence it is sufficient to 
 find $\varphi\in {\bf H}_{\cal R}$ such that $\varphi\kappa 
 \varphi^{-1}(B)$ contains  $ B_i$ for each $\alpha_i\notin\pi_{\Theta}$.  
 
A straightforward computation using (1.2) shows that
$$q^{(\lambda,\alpha_i)}c\tau(\lambda)y_i-y_it_ic\tau(\lambda)t_i^{-1}=[(\adr 
y_i)c]\tau(\lambda)$$ for all $c\in U$, $1\leq i\leq n$, and 
$\tau(\lambda)\in T$.  Thus $[(\adr 
y_i)b\tau(\lambda)^{-1}]\tau(\lambda)t\in B$ for 
any $\tau(\lambda)\in T$, $t\in T_{\Theta}$, $\alpha_i\in\pi_{\Theta}$, and $b\in B.$  
As in [L4, following the proof of Theorem 7.5], we have $$( {\rm ad}_r
y_{i_1}^{(m_1)})\cdots 
({\rm ad}_ry_{i_r}^{(m_r)})
\tilde\theta(y_{p(i)})=(-1)^{(m_1+\cdots +m_r)}t_{i}^{-1}x_{i}$$
and $\kappa((\adr x_j)b)=-((\adr y_j)\kappa(b)$ for all $b\in U$
and $1\leq j\leq n$.  Thus (2.1) implies that     $$\kappa(\tilde\theta(y_i))=
(-1)^{(m_1+\dots+m_r)}(\adr y_{i_1}^{(m_1)})\cdots 
(\adr y_{i_r}^{(m_r)})y_{p(i)}.$$  Moreover, it follows   that 
$$\eqalign{&\kappa(\tilde\theta(y_i))t_i+d(\eta_{p(i)})t_i^{-1}x_it_i+ 
s(\zeta_i)t_i\cr&=(-1)^{(m_1+\cdots+m_r)}[(\adr y_{i_1}^{(m_1)})\cdots 
(\adr 
y_{i_r}^{(m_r)})B_{p(i)}t_{p(i)}^{-1}]t_{p(i)}(t_{p(i)}^{-1}t_i)\cr
}\leqno(3.1)$$
is an element of $B$.

 Let $\varphi$ be the Hopf algebra automorphism in ${\bf 
H}_{\cal R}$ defined by 
$$\varphi(x_j)=q^{(\alpha_j,-\alpha_j+\Theta(\alpha_j))/4} 
d(\eta_j)^{1/2}x_j{\rm\ for\ } 
\alpha_j\notin\pi_{\Theta}.$$ Temporarily write 
$\kappa_1=\varphi\kappa\varphi^{-1}$.
By (2.1) and (2.2), we see that
$\Theta(\alpha_i)-\alpha_i= \Theta(\alpha_{p(i)})-\alpha_{p(i)}$.
Hence 
$$\eqalign{&(\alpha_i,\Theta(\alpha_i)-\alpha_i)=
(\alpha_i,\Theta(\alpha_{p(i)})-\alpha_{p(i)})=
(\alpha_i,\Theta(\alpha_{p(i)}))-(\alpha_i,\alpha_{p(i)})\cr
&=(\Theta(\alpha_i)-\alpha_i,\alpha_{p(i)})
=(\alpha_{p(i)},\Theta(\alpha_{p(i)})-\alpha_{p(i)}).\cr}$$
Thus, a straightforward computation
using the fact that $s(\zeta_i)=0$ for $\alpha_i\notin{\cal S}$
shows that
$$\eqalign{&q^{(\alpha_i,\alpha_i+\Theta(\alpha_i))/2}\kappa_1(\kappa(\tilde\theta(y_i))t_i+
 d(\eta_{p(i)})t_i^{-1}x_it_i+\tilde 
s_it_i)\cr&=d(\eta_{p(i)})  d(\eta_i)^{-1}(y_it_i+  d(\eta_{p(i)})
\tilde\theta(y_i)t_i+ 
s(\zeta_i)t_i)=  d(\eta_{p(i)}) d(\eta_i)^{-1} B_i\cr} $$ is an element of $\kappa_1(B)$ for all
$\alpha_i\in \pi_{\Theta}$.   Therefore
$B=\kappa_1(B)$.
$\Box$ 

\medskip

Recall that $P^+(\Sigma)$ is the subset of $\gh^*$ consisting of dominant integral weights 
associated to the root system $\Sigma$.   Let 
$P^+_{\Theta}$ be the set of all $\lambda\in P^+(\pi)$ 
such that \begin{enumerate}
\item[(i)]$ (\lambda,\beta)=0$ for all $\beta\in \gh^*_{\Theta}$.
\item[(ii)]
$(\tilde\lambda,\tilde\beta)/(\tilde\beta,\tilde\beta)$ is an integer for
every   restricted root $\tilde\beta\in \Sigma$. 
\end{enumerate}
By  [He2, Chapter II, Theorem 4.8], $P^+_{\Theta}$ is just the intersection of $P(2\Sigma)$ with 
$P^+(\Sigma)$.

Given $\lambda\in P^+(\pi)$ and an integral domain ${\cal Z}$ 
containing ${\cal 
C}$, we set 
$L({\lambda})_{{\cal Z}}=L({\lambda})\otimes_{\cal C}{\cal Z}$. We have the following
criterion for $L(\lambda)_{\cal K}$ to be a spherical module with respect to the  
the algebra $B_{{\theta},{\bf s}(\zeta), {\bf d}(\eta)}$.

\begin{theorem}  Let $B= B_{{ \theta},{\bf  
s}(\zeta),{\bf  d}(\eta)}$  Then $$\dim L(\lambda)_{\cal K}^B\leq 
1\leqno{(3.2)}$$ for all $\lambda\in P^+(\pi)$.  Moreover equality holds in (3.2) and 
thus  $L(\lambda)_{\cal K}$    is spherical with respect to 
$B_{{\theta},{\bf  
s}(\zeta),{\bf  d}(\eta)}$   if and only if 
$\lambda\in P^+_{\Theta}$.  
\end{theorem}

\noindent{\bf Proof:}  The first assertion follows as in the proof 
of [L4, Theorem 7.7(i)] and the second assertion follows as in the 
proof of  [L3, Theorem 4.3] using Theorem 3.1. $\Box$

\medskip

    Fix  $\lambda\in 
P^+_{\Theta}$. Given a weight $\mu$ of $L(\lambda)$,  let 
$\{v_{\mu}^i\}_i$ be a basis of weight vectors  for
the $\mu$ weight space of 
$L({\lambda})$ where $i$ varies over a finite set. 
Note that the $\lambda$ weight 
space of $L(\lambda)$ has dimension $1$ and so $\{v_{\lambda}^i\}_i$ 
has only one element.  Set 
$v_{\lambda}=v_{\lambda}^1$.
Given a vector $v=\sum_{i,\mu}g^i_{\mu}v^i_{\mu}$ in $L({\lambda})_{{\cal 
K}}$,
set $${\rm Supp}(v)=\{\mu|g^i_{\mu}\neq 0{\rm \ for \ some \ }i\}.$$

Write ${\cal C}[\zeta,\eta]$ for the 
polynomial ring ${\cal C}[\zeta_i,\eta_j|\alpha_i\in {\cal S}$ 
and $\alpha_j\in {\cal D}]$. 
Set
$$C_k=x_k+q_k^{2}d(\eta_{p(k)})^{-1}\kappa(\tilde\theta(y_k))t_k+q_k^{2} 
s(\zeta_k)(t_k-1)\leqno{(3.3)}$$ for $\alpha_k\notin\pi_{\Theta}$ and
   set  $C_k=x_k$ for 
$\alpha_k\in \pi_{\Theta}$. 
Note that 
$\tilde\theta(y_k)\in U_{{\cal C}[\zeta,\eta]}$ for each 
$\alpha_k\notin \pi_{\Theta}$.  Hence,
by (3.1)
and the definition of the counit of $U$, 
each $C_k$ is an element of $B_+\cap U_{{\cal C}[\zeta,\eta]}$.

The next  result provides  finer information about the invariant 
vectors of spherical modules.

\begin{lemma}
  Suppose that $\tilde\xi_{\lambda}$ 
is a nonzero 
$B_{\theta,{\bf   s}(\zeta),{\bf  d}(\eta)}$ invariant vector of 
$L({\lambda})_{\cal K}$.  Then  
$$\tilde\xi_{\lambda}\in v_{\lambda}+\sum_{\mu<\lambda,k} 
{\cal C}[\zeta,\eta]v_{\mu}^k$$ up to multiplication by  a 
nonzero scalar. 
\end{lemma}

\noindent
{\bf Proof:} Let $a_{\lambda}$ denote the 
coefficient of $v_{\lambda}$ in
$\tilde\xi_{\lambda}$ written as a linear combination of the basis 
vectors $\{v_{\mu}^i\}$ of $L(\lambda)$.  Rescaling if necessary, we may assume that 
$a_{\lambda}$ equals $0$ or $1$.

We can write $\tilde\xi_{\lambda}$ as a 
linear combination of weight vectors:  
$\tilde\xi_{\lambda}= \sum_{\mu}w_{\mu}.$  
For each weight $\gamma$ of $L(\lambda)$, let 
$y_{-\lambda+\gamma}^i$ be weight vectors  in $U^-_{-\lambda+\gamma}$ such 
that $y_{-\lambda+\gamma}^iv_{\lambda}=v_{\gamma}^i.$ It is well known that 
there exists a set of weight vectors $\{x_{\lambda-\gamma}^i\}$   in $U^+_{\lambda-\gamma}$ such that 
$x_{\lambda-\gamma}^iv_{\gamma}^j=\delta_{ij}v_{\lambda}.$
It follows that 
$\sum_iy_{-\lambda+\gamma}^ix_{\lambda-\gamma}^iw_{\gamma}=w_{\gamma}$.   
In particular, $w_{\gamma}\in \sum_iUx_iw_{\gamma}$ for each 
$\gamma\in {\rm Supp}(\tilde\xi_{\lambda})-\{\lambda\}$.

Partition ${\rm Supp}(\tilde \xi_{\lambda})$ into 
two sets $I_1$ and $I_2$ as follows.  The weight $\gamma$ is in    $I_1$ provided 
$w_{\gamma}$ is a nonzero element of 
$a_{\lambda}U_{{\cal C}[\zeta,\eta]}v_{\lambda}$.  The second set $I_2$ is simply the complement of 
$I_1$ in ${\rm Supp}(\tilde\xi_{\lambda})$.   Note that $I_1\cup I_2$ is 
nonempty since $\tilde\xi_{\lambda}$ is nonzero.  We show that $I_2$ is empty.
This in turn implies  that $I_1={\rm Supp}(\tilde\xi_{\lambda})$ and 
$a_{\lambda}=1$, which proves the theorem.

Assume $I_2$ is nonempty and  let $\gamma$ be a maximal element of $I_2$.
Fix $i$. The maximality of the 
choice of $\gamma$ and the definition of $C_k$  ensures that
$$C_k\sum_{\mu\in I_2}w_{\mu}\in 
x_kw_{\gamma}+\sum_{\mu\not\geq\gamma+\alpha_k}(L(\lambda)_{\cal K})_{\mu}.$$ 
Recall that  that $C_k\tilde\xi_{\lambda}=0$.    Thus  
 $$C_k\sum_{\mu\in I_1}w_{\mu}\in 
-x_kw_{\gamma}+\sum_{\mu\not\geq\gamma+\alpha_k}(L(\lambda)_{\cal K})_{\mu}.$$ 
However, since
$C_k\in U_{{\cal C}[\zeta,\eta]}$, it follows that 
$$C_k\sum_{\mu\in I_1}w_{\mu} \in a_{\lambda}U_{{\cal C}[\zeta,\eta]}v_{\lambda}. $$ 
Thus $x_kw_{\gamma}\in a_{\lambda}U_{{\cal C}[\zeta,\eta]}v_{\lambda}. $
for all $1\leq k\leq n$.  It follows that 
$$w_{\gamma}\in \sum_kUx_kw_{\gamma}\in a_{\lambda}U_{{\cal 
C}[\zeta,\eta]}v_{\lambda}, $$ a contradiction.  Therefore $I_2$ is empty.
$\Box$

\medskip
Note that Theorem 3.2 was originally proved in [L3] only for subalgebras  in ${\cal 
B}_{\theta,1}$.   Now consider spherical vectors corresponding to a subalgebra $B$ 
of $U$ 
in the larger set ${\cal B}_{\theta}$. The argument used to sketch the 
proof of   [L4, Theorem 
7.7(i)] shows 
that $$\dim L(\lambda)^B\leq 1$$ for all $\lambda\in P^+(\pi)$.  
Moreover the proof of [L3, Theorem 
4.3] can be applied in this more general setting to show  that if $L({\lambda})$ 
admits a spherical vector then $\lambda\in P^+_{\Theta}$.   However, 
the proof of the existence of a spherical vector associated to $B$ in ${\cal 
B}_{\theta,1}$ relies  on the fact that
the $B\in {\cal B}_{\theta,1}$ act semisimply on 
finite-dimensional $U$ modules.   Unfortunately, we  have not proved a 
version of Theorem 3.1 for all subalgebras  $B$ 
of $U$ 
in ${\cal B}_{\theta}$. Nevertheless, we obtain a version of Theorem 
3.2,  which applies to all the   algebras 
in ${\cal B}_{\theta}$, using the preceding lemma.

\begin{theorem} Let ${\cal Z}$ be an integral domain containing ${\cal C}$
and let  $B$ be an algebra in ${\cal B}_{\theta}({\cal Z})$.
  Then
 $$\dim L(\lambda)_{{\cal Z}}^B\leq 
1.\leqno{(3.4)}$$  Moreover equality holds in (3.4) (and thus  
$L(\lambda)_{{\cal Z}}$    
is spherical with respect to $B$)   if and only if 
$\lambda\in P^+_{\Theta}$.  
\end{theorem}

\noindent{Proof:}   Applying a Hopf 
algebra automorphism of $U_{{\cal Z}}$ to $B$ if necessary, we reduce to the case when 
$B=B_{\theta,{\bf s},{\bf d}}$ for some ${\bf s}=(s_1,\dots, s_n)\in {\bf 
S}({\cal Z})$ 
and ${\bf d}=(d_1,\dots, d_n)\in {\bf D}({\cal Z})$.  Without loss of 
generality, we may assume that ${\cal Z}$ is the ${\cal C}$ algebra 
generated by the $s_i$ for  $\alpha_i\in {\cal S}$ and the $d_j$ for 
$\alpha_j\in {\cal D}$.
 By the discussion preceding 
the theorem, we need only show that there exists a $B$ invariant 
vector in $L({\lambda})_{{\cal Z}}$ for all  
$\lambda\in P^+_{\Theta}$. 

Fix  $\lambda\in 
P^+_{\Theta}$. 
 By Lemma 3.3, there exist polynomials  
$f_{\mu}^k$ in ${\cal C}[\zeta,\eta]$
such that $$\tilde\xi_{\lambda}=v_{\lambda}+\sum_{\mu<\lambda,k} 
f_{\mu}^kv_{\mu}^k$$ is a 
$B_{\theta,{\bf  s}(\zeta),{\bf  d}(\eta)}$ invariant vector. Let $J$ be the 
ideal in ${\cal C}[\zeta,\eta]$ generated 
by $\zeta_i-s_i$ as
$\alpha_i$ runs over ${\cal S}$ and 
$\eta_j-d_j$ as $\alpha_j$ runs over $ {\cal D}$. 
 Let $\bar f_{\mu}^k$ denote 
the image of $f_{\mu}^k$ in ${\cal Z}$  obtained by 
modding out by  the ideal $J$. 
Set $$\xi_{\lambda}=v_{\lambda}+\sum_{\mu<\lambda,k} 
\bar f_{\mu}^kv_{\mu}^k.$$  Note that $\xi_{\lambda}$ is nonzero since the 
coefficient of $v_{\lambda}$ is $1$.
Also, the 
image  
of  $B_{\theta,{\bf  s}(\zeta),{\bf   d}(\eta)}\cap (U_{{\cal 
C}[\zeta,\eta]})_+$ in 
$U_{{\cal Z}}$ obtained by modding out by $J$ is just $B_+$.
Furthermore, the action of  $U_{{\cal Z}}$ on $L(\lambda)_{{\cal Z}}$
induced by taking the action of $U_{{\cal C}[\zeta,\eta]}$ on 
$L(\lambda)_{{\cal C}[\zeta,\eta]}$ and modding 
out by $J$ is the same as the standard action of $U_{{\cal Z}}$ 
on $L(\lambda)_{{\cal Z}}$.
Thus $\xi_{\lambda}$ is a nonzero $B_+$ invariant vector in 
$L(\lambda)_{\cal Z}$.   
$\Box$
  
\medskip

Recall that $P(\Sigma)$ is the weight lattice associated to the root 
system $\Sigma$ and hence contains the root lattice associated to 
$\Sigma$.  
Since 
$\tilde\alpha_i\in \Sigma$ it follows that $\tilde\alpha_i\in 
P(\Sigma)$ for each $\alpha_i\in \pi-\pi_{\Theta}$.
Sometimes, a smaller scalar multiple of $\tilde\alpha_i$ is also in 
$\Sigma$.   This happens when ${\cal S}$ is nonempty as explained in 
the following lemma.

\begin{lemma} If $\alpha_i\in {\cal S}$ then $\tilde\alpha_i/2\in P(\Sigma)$.
\end{lemma}

\noindent
{\bf Proof:} It is sufficient to prove the lemma when $\gg,\theta$ is 
an 
irreducible pair.   Given $\alpha_i\in {\cal S}$, we must check that 
$${{2(\tilde\alpha_i/2,\tilde\alpha_j)}\over(\tilde\alpha_j,\tilde\alpha_j)}\in 
{\bf Z}\leqno{(3.5)}$$  for all $\tilde\alpha_j\in \Sigma.$  This is clearly true 
when $i=j$.  
Hence it is sufficient to check the case when $i\neq j$. 
 By the classification in Section 7,
  there are only  four possibilities with ${\cal S}$ nonempty: {\bf AIII}, Case 2; {\bf CI}; 
{\bf DIII} Case 1;  {\bf EVII}.  
The tables in [A, Section 5] and [He, Chapter X, Section F] include 
the Dynkin diagram associated to the restricted root systems in the 
column labelled ``$\Delta^-$''.   In each of these four cases, the 
restricted root system is of type $C$ and the root $\tilde\alpha_i$, 
for $\alpha_i\in {\cal S}$, corresponds to the unique long root.   
Hence (3.5) follows.
$\Box$

\medskip

Define an ordering on the restricted weights
$\{\tilde\lambda|\lambda\in Q(\pi)\}$ by $\tilde\beta'{\geq}_r\tilde\beta$ 
if and only if $\tilde\beta'-\tilde\beta\in 
\sum_i{\bf N}\tilde\alpha_i$. Now suppose that $\lambda$ and 
$\mu$ are two elements of $P(\pi)$ such that $\lambda\geq \mu$.  Then 
there exist nonegative integers $n_i$ such that 
$\lambda-\mu=\sum_in_i\alpha_i$.  It follows that 
$\tilde\lambda-\tilde\mu=\sum_in_i\tilde\alpha_i$.   Hence 
$\tilde\lambda\geq_r\tilde\mu$ with strict inequality if and only if $n_i$ is 
nonzero for some $\alpha_i\in\pi-\pi_{\Theta}$.

The next result  gives additional information about  the spherical 
vectors which will be used in later sections.

\begin{theorem}
Let ${\cal Z}$ be an integral domain containing ${\cal C}$, let 
$B\in {\cal B}_{\theta}({\cal Z})$ and let 
$\lambda\in P^+_{\Theta}$.  Write  
$\xi_{\lambda}$ for  a nonzero 
$B$ invariant vector of 
$L(\lambda)_{{\cal Z}}$.  Then 
${\rm Supp}(\xi_{\lambda})\subset P(2\Sigma)$,  $\lambda\in {\rm 
Supp}(\xi_{\lambda})$, and 
if $\beta\in {\rm Supp}(\xi_{\lambda})-\{\lambda\}$ then
$\tilde\beta<_r\tilde\lambda$. 
\end{theorem}

\noindent
{\bf Proof:} Without loss of generality, we may assume that 
$B=B_{\theta,{\bf s},{\bf d}}$
for some ${\bf s}\in {\bf S}({\cal Z})$ and ${\bf d}\in {\bf   
D}({\cal Z})$.
By  Theorem 3.4, rescaling if necessary, we can   
write $$\xi_{\lambda}=v_{\lambda}+\sum_{\mu<\lambda} 
w_{\mu}$$ where each $w_{\mu}$ is a weight vector of weight $\mu$ in 
$L(\lambda)_{{\cal Z}}$. 
In particular,
$\lambda\in {\rm Supp}(\xi_{\lambda})$.

Recall that $\{\tilde\alpha_i|\alpha_i\in\pi-\pi_{\Theta}\}$
is a set of simple roots for $\Sigma$.  It follows that 
 that $2\tilde\alpha_i\in P(2\Sigma)$ for all 
$\alpha_i\in \pi-\pi_{\Theta}$. By Lemma 3.5, $\tilde\alpha_i\in 
P(2\Sigma)$ for all $i$ such that $\alpha_i\in{\cal S}$. 
Set $$N=\sum_{\alpha_i\in\pi-\pi_{\Theta}-{\cal S}}
{\bf N}(2\tilde\alpha_i)+\sum_{\alpha_i\in{\cal S}}{\bf 
N}\tilde\alpha_i\leqno{(3.6)}$$
and $$\lambda-N=\{\lambda-\gamma|\gamma\in N\}.$$ Now $\lambda\in P^+_{\Theta}$ and $P^+_{\Theta}$ 
is a subset of $P(2\Sigma)$. Hence $\lambda-N\subset P(2\Sigma)$.
Furthermore, every element 
$\beta\in \lambda-N$ such that $\beta\neq \lambda$ satisfies 
$\beta<_r\lambda$.  We complete the proof of the theorem by
showing   that $${\rm 
Supp}(\xi_{\lambda})\subset\lambda-N.$$

   Consider $\alpha_i\in 
\pi_{\Theta}$ and recall that $x_i\xi_{\lambda}=0$.  Note that if 
$\beta\neq \gamma$, then 
$x_iw_{\beta}$ and $x_iw_{\gamma}$ have the same weight if and only 
if 
they are both zero. 
 Thus 
$x_iw_{\beta}=0$ for each $\beta\in {\rm Supp}(\xi_{\lambda})$
and $\alpha_i\in\pi_{\Theta}$.

Let $ I_1$ be the subset of ${\rm 
Supp}(\xi_{\lambda})$ consisting of those 
weights $\beta$ such that  $\beta\in \lambda-N$.  Since $\lambda\in 
{\rm 
Supp}(\xi_{\lambda})$ and $\lambda\in \lambda-N$, it follows that $\lambda\in I_1$. Let 
$ I_2$ be the complement of $ I_1$ in ${\rm 
Supp}(\xi_{\lambda})$.  We argue that $ I_2$ is 
empty.   Choose $\beta\in I_2$ such that 
$\beta$ is a maximal element in $I_2$ with respect to the partial ordering 
$>_r$.  Since $\lambda\neq\beta$, $w_{\beta}$ is not a highest weight 
vector in $L(\lambda)_{{\cal Z}}$.  Hence there exists $i$ such 
that   
$x_iw_{\beta}\neq 0$. By the previous paragraph, we must have that 
$\alpha_i\in\pi-\pi_{\Theta}$.   Note that $x_iw_{\beta}$ is a weight vector of 
weight $\alpha_i+\beta$. Assume first that $\alpha_i\notin{\cal S}$. It 
follows from the definition of $C_i$  (see (3.3)) that   $C_i=x_i+Y_i$ where  
$Y_i$  is a weight vector in $U$   of weight 
$\Theta(\alpha_i)$.  Now $C_i\xi_{\lambda}=0$ since $C_i\in B_+$.  
Hence  
$x_iw_{\beta}=-Y_iw_{\gamma}$ for some $\gamma\in {\rm 
Supp}(\xi_{\lambda})$.  Moreover 
$\gamma=\alpha_i-\Theta(\alpha_i)+\beta=2\tilde\alpha_i+\beta$.  It follows that 
$\gamma-\beta=2\tilde\alpha_i$ and thus 
$\gamma>_r\beta$. Moreover, $\beta\notin \lambda-N$ implies that
$\gamma\notin \lambda-N$ which contradicts the choice of $\beta$.

Now   assume $\alpha_i\in {\cal S}$. Hence 
$\Theta(\alpha_i)=-\alpha_i$ and $\tilde\alpha_i=\alpha_i$.  Note that $B_+$ 
contains the element $C_i=x_i+q_i^2y_it_i+q_i^2s_i(t_i-1)$.  Hence we can 
write 
$x_iw_{\beta}=-q_i^2y_it_iw_{\gamma_1}-q_i^2s_i(t_i-1)w_{\gamma_2}$ where $\gamma_1=\beta+2\tilde\alpha_i$ 
and $\gamma_2=\beta+\tilde\alpha_i$. Since $x_iw_{\beta}$ is 
nonzero, at least one of $\gamma_1,\gamma_2$ is contained in  ${\rm 
Supp}(\xi_{\lambda})$.  But  neither $\gamma_1$ nor 
$\gamma_2$ can be in the subset $I_1$ of $\lambda-N$ since $\beta\notin\lambda-N$.     Once again this contradicts the choice 
of $\beta$ and the lemma follows.
$\Box$

\medskip

\section {Quantum Zonal Spherical Functions}

In this section, we review some basic facts about the quantized 
function algebra and define quantum zonal spherical functions 
associated to pairs $B,B'$ of algebras in ${\cal B}_{\theta}$.   We then 
prove an important injectivity result relating the space of 
biinvariants
 to the character group ring associated to $\Sigma$.

Let $R_q[G]$ denote the quantized function algebra as defined 
in  [Jo, Chapter 9.1]. Recall that $R_q[G]$ is a right and left $U$  module.  
Let $L(\lambda)^*$ denote the $U$ module dual to $L(\lambda)$.  Viewing $L(\lambda)$ as a left $U$ 
module, $L(\lambda)^*$ is given its 
natural right $U$ module structure.

According to the quantum  Peter-Weyl theorem ([Jo, 9.1.1 and 
1.4.13]), there is 
an isomorphism as right and left $U$ modules:
$$R_q[G]\cong \oplus_{\lambda\in P^+(\pi)} L(\lambda)\otimes 
L(\lambda)^*.\leqno{(4.1)}$$ Given an element $w\otimes w^*\in L(\lambda)\otimes 
L(\lambda)^*$, we write the corresponding element of $R_q[G]$ as 
$c^{\lambda}_{w^*,w}$.  As a vector space, 
$R_q[G]$ is spanned by vectors $c^{\lambda}_{w^*,w}$ where 
$\lambda\in P^+(\pi)$, $w\in 
L(\lambda)$, and $w^*\in L(\lambda)^*$.
Note that elements of $R_q[G]$ can be thought of as 
functions on the quantized enveloping algebra $U$.   In particular, 
$$c^{\lambda}_{w^*,w}(u)=w^*(uw)$$ for all $u\in U$. 

Let  $B$ and $B'$ be two subalgebras of $U$ in ${\cal B}_{\theta}$.
 Define the 
subspace of ${ }_{B'}{\cal H}_{B}$ of left  $B'$ and right $ B$   invariants  
inside of 
 $R_q[G]$ by $$\eqalign{{}_{B'}{\cal H}_{B}
= &\{\varphi\in R_q[G]|b'\cdot \varphi=\epsilon(b')\varphi {\rm\ and \ } 
 \varphi\cdot 
b =\epsilon(b )\varphi\cr & {\rm\ for\ all\ } b\in B{\rm\ and\ }b'\in 
B' \}.\cr}$$
By [KS, Chapter 11, Proposition 68], 
  ${ }_{B'}{\cal H}_{B}$, which we refer to as the space of 
  biinvariants associated to the pair $(B,B')$, is a
subalgebra  of $R_q[G]$.

Set $${ }_{B'}{\cal H}_{B}(\lambda)={  }_{B'}{\cal H}_{B}\cap 
(L(\lambda)\otimes L(\lambda)^*)$$ where we identify  $L(\lambda)\otimes L(\lambda)^*$ with a 
subspace of $R_q[G]$ using (4.1).  
By Theorem 3.4,
${  }_{B'}{\cal H}_{B}(\lambda)$ is $0$ if $\lambda\notin 
P_{\Theta}^+$ and is a one-dimensional trivial left 
$B'$ and right 
$B $ module otherwise.
Thus (4.1) implies 
the  following direct sum decomposition into trivial one-dimensional 
left 
$B'$ and right 
$B $  modules:
$${  }_{B'}{\cal H}_{B}\cong \oplus_{\lambda\in P^+_{\Theta}}\ {  }_{B'}
{\cal H}_{B}(\lambda).\leqno{(4.2)}$$
For each $\lambda\in P^+_{\Theta}$,
nonzero elements of ${ }_{B'}{\cal H}_{B}(\lambda)$ are called zonal 
spherical functions associated to the pair $B,B'$. Fix a nonzero 
vector
$c^{\lambda}_{B,B'}\in {  }_{B'}{\cal H}_{B}(\lambda)$.
It follows from (4.2) that $\{c^{\lambda}_{B,B'}|\lambda\in P^+_{\Theta}\}$ is a 
basis for ${  }_{B'}{\cal H}_{B}$.  

Note that $L(\lambda)\otimes 
L(\lambda)^*$ is the joint eigenspace for the action (either right or 
left) of the center 
of $U$ on $R_q[G]$ with eigenvalue  given by the 
central character of $L(\lambda)$.   Thus the zonal spherical 
functions are just the eigenvectors of ${ }_{B'}{\cal H}_{B}$
 for the action of the center of $U_q(\gg)$.
Moreover, the eigenvalue of $c^{\lambda}_{B,B'}$ with respect to 
the central element $a$ is given by the action of $a$ on 
$L(\lambda)$.

Given $\lambda\in P(\pi)$, define 
the function $z^{\lambda}$ on $T$ by $$z^{\lambda}\cdot 
\tau(\beta)=q^{(\lambda,\beta)}.$$  Write ${\cal C}[P]$ for the
group 
algebra over ${\cal C}$ generated by the multiplicative group $\{z^{\lambda}|\lambda\in 
P(\pi)\}$.   
 Let $\Upsilon$ be the algebra homomorphism  from $R_q[G]$ to 
${\cal C}[P]$ obtained by 
restricting 
elements of $R_q[G]$, considered as  functions on $U$, to  
the subalgebra $U^o$ generated by $T$. 

Suppose that     ${\cal 
Z}$ is an integral domain 
containing ${\cal C}$.  Set ${\cal Z}[P]={\cal C}[P]\otimes_{\cal 
C}{\cal Z}$.   The map $\Upsilon$ extends to a 
homomorphism from 
$R_q[G]\otimes_{\cal C}{\cal Z}$ to ${\cal Z}[P]$
which we also denote by $\Upsilon$. The definitions of  $ 
{}_{B'}{\cal H}_{B}$ and  $c^{\lambda}_{B,B'}$ also extend in the 
obvious way to    pairs of algebras 
$B,B'$ in ${\cal B}_{\theta}({\cal Z})$.  
    
Suppose that $w$ and $w^*$ are weight vectors in 
$L(\lambda),L(\lambda)^*$ respectively.     Assume that the weight of  
$w$ 
 is $\beta$.  For each 
$\tau(\mu)$ in $T$, we have 
$$c^{\lambda}_{w^*,w}(\tau(\mu))=w^*(\tau(\mu)w)=q^{(\beta,\mu)}w^*(w)
=w^*(w)z^{\beta}\cdot\tau(\mu).$$  
In particular, $\Upsilon(c^{\lambda}_{w^*,w})$  is
a scalar multiple  of $z^{\beta}$ where $\beta$ is the weight of $w$.
Note that if the weight of $w$ and $w^*$ differ, then $w^*(w)=0$ and 
hence $\Upsilon(c^{\lambda}_{w^*,w})=0$.

The next lemma transfers information about spherical vectors from the 
last section to the zonal spherical functions.

\begin{lemma}
Let ${\bf s}\in {\bf S}({\cal Z})$ and ${\bf d}\in {\bf D}({\cal Z})$
where ${\cal Z}$ is an integral domain containing ${\cal C}$.  
Suppose that $B=\chi(B_{\theta,{\bf s},{\bf 
d}})$ and $B'=\chi'(B_{\theta,{\bf s}',{\bf d}})$ where $\chi$ and $\chi'$ are 
Hopf algebra automorphisms of $U_{{\cal Z}}$ 
fixing $T$ and ${\cal M}$. 
Then for each $\lambda\in P^+_{\Theta}$, 
$$\Upsilon(c^{\lambda}_{B,B'})=z^{\lambda}+\sum_{\beta<_r\lambda}f_{\beta}z^{\beta}$$ 
up to a nonzero scalar where each
$f_{\beta}\in {\cal Z}$  and  $\beta$ runs over weights of $L(\lambda)$ 
which are contained in 
$  P(2\Sigma)-\{\lambda\}$. 
\end{lemma}

\noindent
{\bf Proof:} Let $\{v^i_{\mu}|\mu\leq \lambda\}$ be a basis 
consisting of weight vectors 
for $L(\lambda)$ as in Section 3.  Let $\{v^{i*}_{\mu}|\mu\leq \lambda\}$ be 
a basis of weight vectors for 
$L(\lambda)^*$ such that  $v_{\mu}^{i*}(v_{\mu}^{j})=\delta_{ij}$.

By Theorem 3.4, $L(\lambda)_{\cal Z}$
contains a $B'$ invarinat vector $\xi$ and $L(\lambda)_{\cal 
Z}^*$ contains a $B$ invariant vector $\xi^*$.
Rescaling if necessary, we may write  $$\xi=v_{\lambda}+ 
\sum_{\beta<_r\lambda}f_{\beta}^iv_{\beta}^i$$ 
and $$\xi^*=v^*_{\lambda}+ 
\sum_{\beta<_r\lambda}f_{\beta}^{i*}v_{\beta}^{i*}$$ where $f_{\beta}^i$ and 
$f_{\beta}^{i*}$ are elements of ${\cal Z}$ for each $\beta$ and $i$. 
Furthermore, by Theorem 3.6, if either $f_{\beta}^{i}$ or $f_{\beta}^{i*}$ is nonzero 
for some $i$, then  $\beta\in 
P(2\Sigma)$.
Set $f_{\beta }=\sum_if_{\beta}^if_{\beta}^{i*}$.  It follows that up to a nonzero scalar, 
$\Upsilon(c^{\lambda}_{B,B'})$ equals 
$z^{\lambda}+\sum_{\beta}f_{\beta}z^{\beta}$ where $\beta$ runs over 
the weights of $L(\lambda)$ 
contained in 
$  P(2\Sigma)-\{\lambda\}$.
$\Box$

\medskip

Consider $\lambda\in 
P^+_{\Theta}$ and $B,B'\in {\cal B}_{\theta}$.  For the remainder of 
the paper, we assume that 
$c^{\lambda}_{B,B'}$ has been chosen so that when 
$\Upsilon(c^{\lambda}_{B,B'})$ is written as a linear combination of the 
$z^{\gamma}$, the coefficient of $z^{\lambda}$ is equal to $1$.
For each $\lambda\in P^+_{\Theta}$, we set
 $\varphi^{\lambda}_{B,B'}=\Upsilon(c^{\lambda}_{B,B'})$.
 
 Set ${\cal C}[2\Sigma]$ equal to the group subalgebra of ${\cal 
 C}[P]$ corresponding to the group $\{z^{2\lambda}|\lambda\in 
P(\Sigma)\}$.  The next result allows us to identify the space 
spanned by the 
zonal spherical functions with a   subalgebra of ${\cal C}[2\Sigma]$

\begin{theorem}
The restriction map $\Upsilon$ from 
${  }_{B'}{\cal H}_{B}$ to  $U^o$ defines an injection of ${  }_{B'}{\cal H}_{B}$
into ${\cal C}[2\Sigma]$.
\end{theorem}

\noindent
{\bf Proof:} Let $\lambda\in P^+_{\Theta}$ and recall that 
$P^+_{\Theta}\subset P(2\Sigma)$. By Lemma 4.1,  there exists 
$a_{\beta}\in {\cal C}$ such that 
$$\varphi_{B,B'}^{\lambda}=z^{\lambda}+\sum_{\beta<_r\lambda}a_{\beta}z^{\beta}$$
where the $\beta$ run over a finite number of elements in $P(2\Sigma)$.
Thus each $\varphi^{\lambda}_{B,B'}\in {\cal C}[2\Sigma]$. It follows 
that $\Upsilon(_{B'}{\cal H}_{B})$ is a subalgebra of ${\cal C}[2\Sigma]$.

Recall that the set $\{c_{B,B'}^{\lambda}|\lambda\in P^+_{\Theta}\}$ 
is a basis for ${  }_{B'}{\cal H}_{B}$.  Consider 
a typical element $X=\sum_ib_ic_{B,B'}^{\lambda_i}$ in  ${\ 
}_{B'}{\cal H}_{B}$ and note that 
$\Upsilon(X)=\sum_ib_i\varphi_{B,B'}^{\lambda_i}$.   
Choose a maximal weight $\lambda_1$ in the set $\{\lambda_i|b_i\neq 
0\}$.    Then $\Upsilon(X)$  is in 
$$b_1z^{\lambda_1}+\sum_{\mu\not\geq_r\lambda_1}{\cal 
C}z^{\mu}$$ and so is nonzero.   Hence $\Upsilon$ is injective.
$\Box$

\medskip

\section{A Criterion for  Invariance}

\medskip   

Let $W$ be the Weyl group associated to the restricted root 
system $\Sigma$. Note that $W$ acts on ${\cal C}[2\Sigma]$ via $w\cdot 
z^{\beta}=z^{w\beta}$ for all $\beta\in P(2\Sigma)$ and $w\in W$.
 A classical result shows that the zonal spherical 
functions correspond to the $W$ invariant functions of the character ring 
of the restricted root system. In this 
section, we obtain a similar result in the quantum case.  In 
particular, given $B\in {\cal B}_{\theta}$ we determine how to choose 
$B'$ so that the image of ${  }_{B'}{\cal H}_{B}$ under $\Upsilon$ is 
the entire invariant ring ${\cal C}[2\Sigma]^{W}$.

Let $\rho$ denote 
the half sum of the positive roots in $\Delta$, so 
$(\rho,\alpha_i)=(\alpha_i,\alpha_i)/2$  for each $1\leq i\leq n$. 
   Note that if  $\Theta(\alpha_i)=-\alpha_{p(i)}$ then
$\tilde\theta(y_i)=t_{p(i)}^{-1}x_{p(i)}$ and $\tilde\theta(y_{p(i)})=t_{i}^{-1}x_{i}$
Hence  $$\tilde\theta(y_{p(i)})t_{p(i)}^{-1}x_{p(i)}=t_{i}^{-1}x_i\tilde\theta(y_i)
 $$ in this case.  We show that a similar 
result is true in general.

Recall that ${\cal M}^+$ is the subalgebra of ${\cal 
M}$ generated by the $x_i$, $\alpha_i\in\pi_{\Theta}$.

\begin{lemma}For each $\alpha_i\notin\pi_{\Theta}$,
$$q^{(\rho,\Theta(\alpha_i)+\alpha_i)}
\tilde\theta(y_{p(i)})t_{p(i)}^{-1}x_{p(i)}\in t_i^{-1}x_i\tilde\theta(y_i)
+ {\cal 
M}^+_+U+U{\cal M}^+_+.$$
\end{lemma}

\noindent
{\bf Proof:}
Recall the sequences used in Section 2 (see (2.1) and (2.2))   to  define 
$\tilde\theta(y_i)$ for $\alpha_i\notin\pi_{\Theta}$. 
Using the form of $(\adr x_i)$ given in (1.2), we have
that $$\eqalign{&t_i^{-1}x_i\tilde\theta(y_i)=t_i^{-1}x_i 
[(\adr x_{i_1})^{(m_1)}\cdots 
(\adr x_{i_r})^{(m_r)}t_{p(i)}^{-1}x_{p(i)}]\cr &\in\
(-1)^{(m_1+\cdots +m_r)}t_i^{-1}x_i [(t_{i_1}^{-1}x_{i_1})^{(m_1)}\cdots 
(t_{i_r}^{-1}x_{i_r})^{(m_r)}t_{p(i)}^{-1}x_{p(i)}]+U{\cal 
M}^+_+.\cr}\leqno{(5.1)}$$ A straightforward computation shows that 
for each positive integer $k$ and each $j$, 
$$(t_{j}^{-1}x_j)^k=q^{-k(k+1)(\alpha_j,\alpha_j)/2}x_j^kt_j^{-k}.\leqno{(5.2)}$$
Note that $x_{i_{s+1}}^{m_{s+1}}\cdots 
x_{i_r}^{m_r}x_{p(i)}$ and 
$X_{s+1}$ (defined in (2.3)) have the same weight.
Recall that $X_s$ is a highest weight vector for the action of $\adr 
x_{i_s}$ and that $X_{s+1}$ is a lowest weight vector for the action 
of $\adr y_{i_s}$.  Set $\lambda_s$ equal to the weight of $X_s$. It follows that 
$(\lambda_s,\alpha_{i_{s}})={{m_s}\over{2}}(\alpha_{i_s},\alpha_{i_s})$ 
and 
$(\lambda_{s+1},\alpha_{i_{s}})=-{{m_s}\over{2}}(\alpha_{i_s},\alpha_{i_s})$
for each $s$.  Hence (5.2) implies 
$$\eqalign{(t_{i_s}^{-1}x_{i_s})^{m_s}x_{i_{s+1}}^{m_{s+1}}&\cdots 
x_{i_r}^{m_r}x_{p(i)}\cr =&q^{-m_s(m_s+1)(\alpha_{i_s},\alpha_{i_s})/2}
x_{i_s}^{m_s}t_{i_s}^{-m_s}
x_{i_{s+1}}^{m_{s+1}}\cdots 
x_{i_r}^{m_r}x_{p(i)}\cr
=&q^{-m_s(\rho,
\alpha_{i_s})}x_{i_s}^{m_s}\cdots 
x_{i_r}^{m_r}x_{p(i)}t_{i_s}^{-m_s}.\cr}\leqno{(5.3)}$$

Note that $$-(\sum_{j=1}^{r} m_j\alpha_{i_j})=\Theta(\alpha_i)-\alpha_{p(i)}.$$
Hence  repeated applications of (5.3) yields
$$\eqalign{(t_{i_1}^{-1}x_{i_1})^{(m_1)}&\cdots 
((t_{i_r}^{-1}x_{i_r})^{(m_r)}t_{p(i)}^{-1}x_{p(i)}\cr&=
 q^{(\rho,\Theta(\alpha_i)-\alpha_{p(i)})}
 x_{i_1}^{(m_1)}\cdots 
x_{i_r}^{(m_r)}x_{p(i)}t_{i_1}^{-m_1}\cdots 
t_{i_r}^{-m_r}t_{p(i)}^{-1}\cr}\leqno(5.4)$$

Now 
$$\eqalign{(\Theta(\alpha_i)+\alpha_{p(i)},\Theta(\alpha_i)-\alpha_i)&=
(\Theta(\alpha_i),\Theta(\alpha_i)-\alpha_i+\alpha_{p(i)})-(\alpha_{p(i)},\alpha_i)\cr
&=(\Theta(\alpha_i),\Theta(\alpha_{p(i)}))-(\alpha_{p(i)},\alpha_i)=0.\cr}$$ Furthermore, 
$$(\rho,\Theta(\alpha_i)-\alpha_{p(i)})+(\alpha_i,\alpha_i)=(\rho,\Theta(\alpha_i)+\alpha_i).$$
Set $m=m_1+\dots+m_r$.  Expressions (5.1) and (5.4) 
imply that $t_i^{-1}x_i\tilde\theta(y_i)$ is an 
element of 
$$
(-1)^{m}q^{(\rho,\Theta(\alpha_i)+\alpha_i)}t_{i}^{-1}t_{i_r}^{-m_r}\cdots 
t_{i_1}^{-m_1}x_ix_{i_r}^{(m_r)}\cdots 
x_{i_1}^{(m_1)}t_{p(i)}^{-1}x_{p(i)}+U{\cal M}^+_+.$$
Using (1.2)  we obtain
$t_{i}^{-1}t_{i_r}^{-m_r}\cdots 
t_{i_1}^{-m_1}x_ix_{i_r}^{(m_r)}\cdots 
x_{i_1}^{(m_1)}t_{p(i)}^{-1}x_{p(i)}$ is an element of 
$$[(\adr x_{i_r}^{(m_r)})
\cdots 
(\ad_rx_{i_1}^{(m_1)})t_i^{-1}x_i]t_{p(i)}^{-1}x_{p(i)}
+{\cal M}^+_+U$$  Thus
$$t_i^{-1}x_i\tilde\theta(y_i)\in 
q^{(\rho,\Theta(\alpha_i)+\alpha_i)}\tilde\theta(y_{p(i)})
t_{p(i)}^{-1}x_{p(i)}+{\cal M}^+_+U+U{\cal M}_+^+.$$
$\Box$

\medskip
The next result provides a criterion for determining when 
${  }_{B'}{\cal H}_{ B}$ is $W$ 
invariant.  In order to prove the result, we need to extend $U$ to a larger algebra.
First, enlarge  $T$ to  the group $\check T$ generated by
$t_i^{1/m}$ for each $1\leq i\leq n$ and each positive integer $m$. 
We can extend the isomorphism $\tau$ to an isomorphism between 
$\sum_i{\bf Q}\alpha_i$ and $\check T$.
 Let 
$\check U$ be the ${\cal C}$ algebra generated by $U$ and $\check T$  such that 
$$\tau(\lambda)v_{\mu}=q^{(\mu,\lambda)}v_{\mu}\tau(\lambda)$$ for 
each $\tau(\lambda)\in \check T$ and weight vector $v_{\mu}\in U$ of 
weight $\mu$.  
    
 Note that each finite-dimensional simple $U$-module $L(\lambda)$ easily becomes 
a $\check U$-module as follows.  Recall that $L(\lambda)$ can be 
written as a direct sum of weight spaces. 
  Let 
$v$ be a weight vector in $L(\lambda)$ of weight $\beta$.  
 Set $$\tau(\mu)\cdot 
v=q^{(\beta,\mu)}v$$ for all $\tau(\mu)\in \check T$. We can 
similarly extend the action of $U$  on $L(\lambda)^*$ to $\check U$. 
It follows that 
elements of $R_q[G]$ extend 
to functions on $\check U$.    
Given an algebra $B$ in ${\cal B}_{\theta}$, let $\check B$ 
 be  the ${\cal C}$ subalgebra of $\check U$ generated by 
$B$ and the group $\check T_{\Theta}=
\{\tau(\lambda)\in\check T|\Theta(\lambda)=\lambda\}.$ 
Note that ${  }_{B'}{\cal H}_{ B}$ is invariant as a left $\check B'$ 
and right $\check B$ module.

   Note that the subalgebra $\check U^o$ of $\check U$ generated by 
   elements of $\check T$
is a  ${\cal C}[P]$ module with the following action: 
$$\sum c_{\gamma}z^{\gamma}\cdot \tau(\lambda)=\sum 
c_{\gamma}q^{(\gamma,\lambda)}.$$ 

\begin{lemma}
Fix $j$ such that $\alpha_j\notin\pi_{\Theta}$ and let 
$\tilde s_{j}$ denote the reflection in $W$ 
corresponding to $\tilde\alpha_j$. Assume that for all  
 $\lambda\in \sum_i{\bf Q}\tilde\alpha_i$ such that
 $(\lambda,\tilde\alpha_j)=0$   that
 $$\tau(\lambda)(\tau(k\tilde\alpha_j)-\tau(-k\tilde\alpha_j))\in \check 
 B_+  \check U+  \check U
\check B'_+\leqno(5.5)$$
for some nonzero rational number $k$.
Then $\Upsilon(_{B'}{\cal H}_{ B})$ is $\tilde s_i$ invariant.
\end{lemma}

\noindent{\bf Proof:} 
Recall that $Q(\Sigma)$ is the root lattice associated to $\Sigma$.
Suppose that  $\beta\in 
 Q(\Sigma)$ and $(\beta,\tilde\alpha_j)\neq 0$.    
 Set $r=(\beta,\tilde\alpha_j)/(\tilde\alpha_j,\tilde\alpha_j)$ and 
 $\lambda=k\beta/r-k\tilde\alpha_j$.  Then
$\beta=(r/k)(\lambda+k\tilde\alpha_j)$, $\lambda\in 
 \sum_i{\bf Q}\tilde\alpha_i$, and $(\lambda,\tilde\alpha_j)=0$.  Hence
$\tilde s_j(\beta)=(r/k)(\lambda-k\tilde\alpha_j)$.

  By (5.5), 
$\tau(k\beta/r)-\tau((k\tilde s_j\beta)/r)\in 
\check B_+\check U+\check U\check B'_+$.  Let $X$ be an element 
of   ${  }_{B'}{\cal H}_{  B}$ and set 
$\Upsilon(X)=\sum_{\gamma}b_{\gamma}z^{\gamma}$.
Note that $X(a)=\Upsilon(X)\cdot(a)$ for all $a\in \check T$ and $X(b)=0$ 
for all $b\in \check B_+  \check U+ \check U
\check B'_+$.   It follows that 
$$\sum_{\gamma}b_{\gamma}z^{\gamma}\cdot 
(\tau(k\beta/r)-\tau((k\tilde s_j\beta)/r))=0.$$
Hence $$\sum_{\gamma}b_{\gamma}z^{({\gamma},k\beta/r)}=
\sum_{\gamma}b_{\gamma}z^{({\gamma},(k\tilde s_j\beta)/r)}.$$
This forces $$\sum_{\{\gamma|(\gamma,k\beta/r)=s\}}b_{\gamma}=
\sum_{\{\gamma|(\gamma,k\tilde s_j\beta)/r)=s\}}b_{\gamma}\leqno(5.6)$$ for each 
integer $s$.  We can rewrite (5.6) as 
$$\sum_{\{\gamma|(\gamma,\beta)=m\}}b_{\gamma}=
\sum_{\{\gamma|(\gamma,\tilde s_j\beta))=m\}}b_{\gamma}$$ for each 
integer $m$. Thus we may conclude that 
$$\sum_{\gamma}b_{\gamma}z^{({\gamma},\beta)}=
\sum_{\gamma}b_{\gamma}z^{({\gamma},\tilde s_j\beta)}=
\sum_{\gamma}b_{\gamma}z^{(\tilde s_j{\gamma},\beta)}\leqno(5.7)$$  for all $\beta$
such that $\beta\in Q(\Sigma)$ and $(\beta,\tilde\alpha_j)\neq 0$. 
Now suppose that $\beta\in Q(\Sigma)$ and  $(\beta,\tilde\alpha_j)= 0$. 
 Then $\tilde s_j\beta=\beta$ and so (5.7) also holds in this case. 
 Since $\sum_{\gamma}b_{\gamma}z^{\gamma}$ is in ${\cal C}[2\Sigma]$, 
 it follows  
that $\sum_{\gamma}b_{\gamma}z^{\gamma}$ is $\tilde s_j$ 
invariant. 
$\Box$

\medskip
Given $B$ in ${\cal B}_{\theta}$, we show how to choose 
another algebra $B'$ in ${\cal B}_{\theta}$ such that   
${  }_{B'}{\cal H}_{ B}$ is $W$ 
invariant.
For each ${\bf c}\in {\bf D}({\cal C})$, 
define the Hopf algebra automorphism $\chi_{{\bf c}}$ 
in ${\bf H}$  of $U$  as follows. 
For all $i$ such that $\alpha_i\in\pi-\pi_{\Theta}$, set
$$\chi_{{\bf c}}(x_i)=q^{-1/2(\rho,\theta(\alpha_i)-\alpha_i)} c_{i}^{-1}x_i$$
and
$$\chi_{{\bf c}}(y_i)=q^{1/2(\rho,\theta(\alpha_i)-\alpha_i)}
c_{i}y_i.$$

Given ${\bf a}$ and ${\bf b}$ in ${\bf D}({\cal C})$ and integers $m,r$, set ${\bf 
a}^m{\bf b}^r$ equal to the $n$-tuple in ${\bf D}({\cal C})$ with entries $a_i^mb_i^r$.

\begin{theorem}
Let $B=B_{{\theta},{\bf s},{\bf d}}$ be in $ {\cal 
B}_{\theta}$. If $B'=\chi_{{\bf c}}(B_{{\theta},{\bf s}',{\bf c}^{2}{\bf d}})$ 
for some ${\bf c}\in {\bf D}({\cal C})$ and ${\bf s}'\in {\bf 
S}({\cal C})$, then  
$\Upsilon(_{ B'}{\cal H}_{B})$ 
  is $W$ invariant.   
\end{theorem}

\noindent
{\bf Proof:} Recall the polynomial ring ${\cal C}[\zeta,\eta]$, its
algebraic closure  ${\cal K}$, and the elements  ${\bf s}(\zeta)\in {\bf 
S}({\cal K} )$ and  $ {\bf 
d}(\eta)\in {\bf D}({\cal K})$  defined in Section 3.
Given ${\bf c}\in {\bf D}({\cal C})$, note that the ${\cal C}$ Hopf algebra automorphism
 $\chi_{{\bf c}}$ of $U$ extends to a ${\cal K}$ Hopf algebra automorphism of 
$U_{\cal K}$ which we also refer to as $\chi_{{\bf c}}$.  Choose ${\bf c}\in {\bf 
 D}({\cal C})$ 
 and ${\bf s}'\in {\bf S}({\cal C})$. 
Set 
 $B=B_{\theta, {\bf s}(\zeta), {\bf d}(\eta)}$ and $B'= 
 \chi_{{\bf c}}(B_{{\theta},{\bf s}',{\bf c}^{2} {\bf d}(\eta)})$.
 By Lemma 4.1 and the proof of Theorem 3.4,  it is enough 
to show that 
$\Upsilon(_{ B'}{\cal H}_{B})$   is $W$ invariant.      

Recall the definition of $\pi^*$ from Section 2.
Fix $j$ such that $\alpha_j\in\pi^*$.   Let $\lambda$ be  an element 
in $\sum_i{\bf Q}\tilde\alpha_i$ such that 
 $(\lambda,\alpha_j)=0$. Since $\tilde\lambda=\lambda$, it follows that  
 $\Theta(\lambda)=-\lambda$ and thus
 $(\lambda,\Theta(\alpha_j))=0$. Thus 
 $\tau(\lambda)$ commutes with both $y_j$ and 
 $\tilde\theta(y_j)$.

 Assume first that 
${\cal S}$ is empty.
Note that
$$B_k=y_kt_k+d(\eta_k)\tilde\theta(y_k)t_k$$
is in $B$ for $k=j$ and $k=p(j)$. Recall (see proof of Theorem 3.1)
that 
$ \Theta(\alpha_j)-\alpha_j=\Theta(\alpha_{p(j)})-\alpha_{p(j)}$. Since $j\leq p(j)$, both 
$c_{p(j)}$ and $d(\eta_{p(j)})$ are equal to $1$. Hence
it follows from the definition of 
$\chi_{{\bf 
c}}$ and (2.1) that $\chi_{\bf 
c}(\tilde\theta(y_j))=q^{-1/2(\rho(\Theta(\alpha_{j})-\alpha_{j})}\tilde\theta(y_j)$.
Moreover,
$$B_j'=(c_{j}
y_jt_j+q^{-(\rho,\Theta(\alpha_j)-\alpha_j)}c_j^{2}d(\eta_j)\tilde\theta(y_{j})t_{j})$$
and $$B_{p(j)}'=(
y_{p(j)}t_{p(j)}+q^{-(\rho,\Theta(\alpha_j)-\alpha_j)}
c_j^{-1}\tilde\theta(y_{p(j)})t_{p(j)})$$
are elements of  $B'$.

By Lemma 5.1, we have
$$\eqalign{&B_kt_k^{-1}x_k=(y_kt_k+d(\eta_k)\tilde\theta(y_k)t_k)t_k^{-1}x_k\cr
\in &{{-(t_k-t_k^{-1})}\over{(q_k-q_k^{-1})}}+t_k^{-1}x_ky_kt_k+ 
q^{-(\rho,\Theta(\alpha_k)-\alpha_k)
}d(\eta_k)t_{p(k)}^{-1}x_{p(k)}\tilde\theta(y_{p(k)})t_{k}
\cr &+{\cal M}^+_+U+U{\cal 
M}^+_+\cr}\leqno{(5.8)}$$ for $k=j$ and $k=p(j)$.   

Consider first the case when $j=p(j)$. In particular,
  $\alpha_j\notin{\cal D}$ 
and thus $c_j=d(\eta_j)=1$. It follows that
$$B_j'=
y_jt_j+q^{-(\rho,\Theta(\alpha_j)-\alpha_j)}\tilde\theta(y_j)t_j.$$
Furthermore,  by (5.8), $$B_jt_j^{-1}x_j\in 
{{-(t_j-t_j^{-1})}\over{(q_j-q_j^{-1})}}+t_j^{-1}x_jB_j'
+{\cal M}^+_+U+U{\cal 
M}^+_+.$$
 Since $\tau(\lambda)$ commutes with $y_j$
 and $\theta(y_j)$, we obtain
$$\eqalign{
&(B_jt_j^{-1}x_j)
\tau(\lambda)-\tau(\lambda)(
t_j^{-1}x_jB_{j}')
\cr&\in 
({{-2(t_j-t_j^{-1})}\over{(q_j-q_j^{-1})}})\tau(\lambda)
+{\cal 
M}^+_+U+U{\cal 
M}^+_+.\cr}$$ Note that 
$(t_j^{-1/2}\tau(\Theta(\alpha_j))^{-1/2}-1)\in \check B_+$ since 
$t_j^{-1/2}\tau(\Theta(\alpha_j))^{-1/2}\in \check T_{\Theta}$.
Therefore
$$\eqalign{&\tau(\lambda)(t_j^{1/2}\tau(\Theta(\alpha_j))^{-1/2}
-t_j^{-1/2}\tau(\Theta(\alpha_j))^{1/2})\cr&=\tau(\lambda)((t_j-t_j^{-1})+
(t_j^{-1/2}\tau(\Theta(\alpha_j))^{-1/2}-1)(t_j+t_j^{-1/2}\tau(\Theta(\alpha_j
)^{1/2}))
\cr&\in \check B_+\check U+\check U\check B'_+.\cr}\leqno{(5.9)}$$
Thus by Lemma 5.2, ${  }_{ B'}{\cal H}_{B}$ is $\tilde s_j$ invariant.

Now consider the case when $j<p(j)$. Note that $t_kt_{p(k)}^{-1}-1$ 
and $t_k^{-1}t_{p(k)}-1$ are both elements of $({\cal 
M}T_{\Theta})_+$ for each $1\leq k\leq n$.   Hence $(t_{p(j)}-t_{p(j)}^{-1})t_{p(j)}^{-1}t_j=(t_j-t_j^{-1})+ ({\cal 
M}T_{\Theta})_+$. Set $a=c_j^{-1}d(\eta_j)^{-1}$ and note that $a+1$ is 
nonzero since $d(\eta_j)$ is an indeterminate. Using (5.8), a straightforward 
computation yields 
$$\eqalign{&aB_jt_j^{-1}x_j+B_{p(j)}t_{p(j)}^{-1}x_{p(j)}t_{p(j)}^{-1}t_j
+(a+1){{(t_j-t_j^{-1})}\over{(q_j-q_j^{-1})}}\cr
&\in c_j^{-1}at_j^{-1}x_jB_j'
+t_{p(j)}^{-1}x_{p(j)}B'_{p(j)}t_{p(j)}^{-1}t_j
+({\cal M}T_{\Theta})^+_+U+U({\cal 
M}T_{\Theta})^+_+.\cr}$$
In particular, (5.9) holds in this case as well and the image of  ${  }_{ B'}{\cal H}_{B}$ is $s_{\tilde\alpha_j}$ invariant
when $j$ is strictly less than $ p(j)$.

Now assume that $\Theta(\alpha_i)=-\alpha_i$ and $\alpha_i\in{\cal 
S}$. 
Then $B_+$ contains $$C_i=y_it_i+q_i^{-2}x_i+  s(\zeta_i)(t_i-1)$$ and
$B'_+$ contains $$C'_i=y_it_i+x_i+s'_i(t_i-1).$$ Furthermore
$$(  s(\zeta_i)-q_i^{-1}s_i')(t_i^{1/2}-t_i^{-1/2})\tau(\lambda)=
C_it_i^{-1/2}\tau(\lambda)-q_i^{-1}\tau(\lambda)
t_i^{-1/2}C_i'.$$ Note that $(
s(\zeta_i)-q_i^{-1}s_i')$ is nonzero since $s(\zeta_i)$ is an indeterminate 
while $s_i'\in {\cal C}$.  It follows that
$\tau(\lambda)(t_i^{1/2}-t_i^{-1/2})\in 
\check B_+\check U+\check U\check B_+'$.  Another application of 
Lemma 5.2 yields that  the image of ${  }_{ B'}{\cal H}_{B}$ under 
$\Upsilon$ is 
$\tilde s_j$ invariant in this case.

We have shown that the image of ${  }_{ B'}{\cal H}_{B}$ under 
$\Upsilon$
is $\tilde s_j$ invariant for all reflections $\tilde s_j$ in the set  
$\{\tilde s_j|\alpha_j\in\pi^*\}$. 
Note that the set $\{\tilde\alpha_j|\alpha_j\in\pi^*\}$ is equal to the set of simple roots of $\Sigma$.  The theorem 
now follows 
from the fact $W$ is generated by the reflections 
in $\{\tilde s_j|\alpha_j\in\pi^*\}$.
$\Box$

\medskip
For each $\lambda\in P^+_{\Theta}$, set $$m_{\lambda}=\sum_{w\in 
W}z^{2w\lambda}.$$ 
It is well known that the set $$\{m_{\lambda}|\lambda\in 
P^+_{\Theta}\}$$
forms a basis for ${\cal C}[2\Sigma]^{W}$.
Now consider a   subset   
$\{\varphi_{\lambda}|\lambda\in P^+_{\Theta}\}$ of ${\cal C}[2\Sigma]$ which satisfies the following 
conditions. For each $\lambda\in P^+_{\Theta}$
\begin{enumerate}
\item[(5.10)] 
$\varphi_{\lambda}=m_{\lambda}+\sum_{\mu<_r\lambda}a_{\lambda,\mu}m_{\mu}$
for some $a_{\lambda,\mu}\in {\cal C}$.
\item[(5.11)]
$\varphi_{\lambda}\in {\cal C}[2\Sigma]^{W}$
\end{enumerate}
Then the set $\{\varphi_{\lambda}|\lambda\in P^+_{\Theta}\}$ forms a basis 
for ${\cal C}[2\Sigma]^{W}$. 

Suppose that $B,B'\in {\cal B}_{\theta}$ are chosen as in Theorem 5.3.
The next result shows that the image of the zonal spherical functions 
associated to the pair $B,B'$ in ${\cal C}[2\Sigma]$ satisfies (5.10)
and  (5.11).

\begin{corollary} Let $B=B_{\theta,{\bf s},{\bf d}}$ be in 
${\cal B}_{\theta}$ and suppose that $B'=\chi_{\bf 
c}(B_{\theta,{\bf s}', {\bf c}^2{\bf d})})$ for some ${\bf 
c}\in {\bf D}$ and ${\bf s}'\in{\bf D}$.  Then the set 
$\{\varphi^{\lambda}_{B,B'}|\lambda\in P^+_{\Theta}\}$ satisfies 
condition (5.10) and (5.11).  Moreover, $\Upsilon$ 
defines an isomorphism of ${  }_{ B'}{\cal H}_{B}$ 
 onto ${\cal C}[2\Sigma]^{W}$.  
\end{corollary}

\noindent
{\bf Proof:} 
By Theorem 4.2, $\varphi^{\lambda}_{B,B'}$ is an element of 
${\cal C}[2\Sigma]$ for each $\lambda\in P^+_{\Theta}$.  On the other hand, 
Theorem 
5.3 implies that each $\varphi^{\lambda}_{B,B'}$ is $W$ 
invariant and hence is an element of ${\cal C}[2\Sigma]^{W}$. It follows that the set $\{\varphi^{\lambda}_{B,B'}|\lambda\in P^+_{\Theta}\}$ satisfies 
condition (5.11).  

Fix $\lambda\in P^+_{\Theta}$. Since $\varphi^{\lambda}_{B,B'}$ is in 
${\cal C}[2\Sigma]^{W}$,  
we can  write $\varphi^{\lambda}_{B,B'}$ as a linear combination 
of the $m_{\gamma},\gamma\in P^+_{\Theta}$:
$$\varphi^{\lambda}_{B,B'}=\sum_{\gamma\in 
P^+_{\Theta}}r_{\gamma}m_{\gamma}.$$ By Lemma 4.1, we can also write 
$$\varphi^{\lambda}_{B,B'}=z^{ \lambda}+
\sum_{ \beta<_r \lambda}c_{ \beta}z^{ \beta}\leqno(5.12)$$ 
for some choice of scalars $c_{ \beta}$ where $\beta$ runs over 
weights of $L(\lambda)$ contained in $ P(2\Sigma)-\{\lambda\}$. Hence, 
if $r_{\gamma}\neq 0$, then
$\gamma$ is a weight of $L(\lambda)$ contained in $ P^+_{\Theta}$. Thus $r_{\gamma}\neq 0$ 
implies that $\gamma=\lambda$ or $\gamma<_r\lambda$. Since $z^{\lambda}$ 
appears with coefficient $1$ in (5.12), it follows that $r_{\lambda}=1$.
Therefore,
 $\varphi^{\lambda}_{B,B'}$  satisfies 
condition (5.10). The last assertion  now follows 
from the discussion preceding the corollary and Theorem 4.2. 
 $\Box$   
\medskip

\section{Zonal Spherical Families}

We say that a function $\lambda\rightarrow \varphi_{\lambda}$ from 
$P^+_{\Theta}$ to ${\cal C}[2\Sigma]$ is a zonal spherical 
family  associated to ${\cal B}_{\theta}$
if there exists $B$
and $B'$ in ${\cal B}_{\theta}$ such that 
 $\varphi_{\lambda}=\varphi^{\lambda}_{B,B'}$ for all $\lambda\in 
 P^+_{\Theta}$.  For may practical purposes, we may think of a zonal 
 spherical family as a set $\{\varphi_{\lambda}|\lambda\in 
 P^+_{\Theta}\}$ indexed by $P^+_{\Theta}$.
Let ${\cal F}_{\theta}$ denote the set of zonal spherical families 
associated to ${\cal B}_{\theta}$. The group ${\bf H}\times {\bf H}$ 
acts on ${\cal B}_{\theta}\times {\cal B}_{\theta}$ in the obvious 
way:
$$(\chi_1,\chi_2)\cdot (B,B')=(\chi_1(B),\chi_2(B')).$$
In this section, we study the relationship 
 between orbits of ${\cal B}_{\theta}\times {\cal B}_{\theta}$ under 
 the action of  ${\bf H}\times {\bf H}$ and orbits of ${\cal F}_{\theta}$ 
 under a subgroup of ${\rm Hom}(Q(\pi),{\cal C}^{\times})$.   This allows 
 us to attach a unique $W$ invariant zonal spherical 
 family (up to possible sign) to each orbit of ${\cal B}_{\theta}\times 
 {\cal B}_{\theta}$.

Recall the definition of the lattice 
$N$ generated by $2\Sigma$ and the set 
$\{\tilde\alpha_i|\alpha_i\in {\cal S}\}$ given in (3.6).
Consider a Hopf algebra automorphism $\chi$ of $U$ in ${\bf H}$. Note 
that
 $\chi$ induces a group homomorphism ${\bar\chi}$ from $Q(\pi)$ to the 
nonzero elements of ${\cal C}$ such that 
$${\chi}(u)= {\bar\chi}(\wt u)u$$ for each weight vector $u\in U$. 
  Since
$N$ is a sublattice of $Q(\pi)$,  $\bar\chi$ restricts to an 
element of ${\rm 
Hom}(N, {\cal C}^{\times})$.

\begin{lemma} The map  $\Psi:{\bf H}\times {\bf H}\rightarrow{\rm 
Hom}(N, {\cal C}^{\times})$ defined by 
$(\chi_1,\chi_2)\rightarrow \bar\chi_1^{-1}\bar\chi_2$ is a surjective
homomorphism.
\end{lemma}

\noindent
{\bf Proof:}  Note first that ${\bf H}$ is an abelian group.   Hence
the map $(\chi_1,\chi_2)\rightarrow  \chi_1^{-1} \chi_2$ defines a 
group  homomorphism 
of ${\bf H}\times {\bf H}$ onto ${\bf H}$.  Thus it is sufficient to 
show that the map $\chi\rightarrow \bar\chi$ defines a surjective 
homomorphism from ${\bf H}$ onto ${\rm 
Hom}(N, {\cal C}^{\times})$. 

Note that $N$ is just the free abelian group with generating set 
$(\pi-\pi_{\Theta}-{\cal S})\cup\{\tilde\alpha_i/2|\alpha_i\in {\cal 
S}\}$.   Set ${{1}\over{2}}N=\{{{1}\over{2}}\beta|\beta\in N\}$.
Let $\psi$ be an element of ${\rm 
Hom}(N, {\cal C}^{\times})$. Taking squareroots allows us to extend 
$\psi$ to a an 
element of  ${\rm 
Hom}({{1}\over{2}}N, {\cal C}^{\times})$.  Set $Q({{1}\over 
{2}}\pi)^{\Theta}$ equal to the set $\{{{1}\over {2}}\beta|\beta\in 
Q(\pi)$ and $\Theta(\beta)=\beta\}$.  We may further extend $\psi$ 
to an element $\hat\psi $ of ${\rm 
Hom}({{1}\over{2}}N+Q({{1}\over 
{2}}\pi)^{\Theta}, {\cal C}^{\times})$ by insisting that 
$\hat\psi(\gamma)=0$ for all $\gamma\in Q({{1}\over 
{2}}\pi)^{\Theta}$. Since $Q(\pi)\subset {{1}\over{2}}N+Q({{1}\over 
{2}}\pi)^{\Theta}$, it follows that $\hat\psi$ restricts to an 
element of ${\rm 
Hom}(Q(\pi) , {\cal C}^{\times})$ which is zero on $Q(\pi)^{\Theta}$.   In 
particular, $\hat\psi=\bar\chi$ where $\chi\in {\bf H}$ 
is chosen so that $\chi(u)=\hat\psi({\rm wt} u)u$ for all weight 
vectors $u\in U$.    $\Box$
\medskip

Note that the action of ${\bf H}\times {\bf H}$ on ${\cal 
B}_{\theta}\times {\cal B}_{\theta}$ induces an action of ${\bf 
H}\times {\bf H}$ on ${\cal F}_{\theta}$.   In particular, given $B$ 
and $B'$ in ${\cal B}_{\theta}$ and $(\chi_1,\chi_2)\in {\bf H}\times 
{\bf H}$, set 
$$(\chi_1,\chi_2)\cdot  \varphi^{\lambda}_{B,B'} = 
\varphi^{\lambda}_{\chi_1(B),\chi_2(B')}$$ for all $\lambda\in 
P^+_{\Theta}$

The group ${\rm Hom}(P(\Sigma),{\cal C}^{\times})$ acts on ${\cal 
C}[2\Sigma]$ via  $g\cdot z^{\beta}=g(\beta)z^{\beta}$ for all $g\in {\rm 
Hom}(P(\Sigma),{\cal 
 C}^{\times})$ and $\beta\in P(2\Sigma).$ This restricts to an action of the 
 group $ {\rm Hom}(N,{\cal C}^{\times})$ on the group ring
${\cal C}[N]$  generated by 
$z^{\beta},\beta\in N$. 

 Let 
$\{\varphi_{\lambda}|\lambda\in P^+_{\Theta}\}$  be a zonal spherical family in ${\cal F}_{\theta}$.  By 
Theorem 3.6 and its proof,  $\varphi_{\lambda}$ is the unique element of the 
vector space ${\cal C}\varphi_{\lambda}$ such that $$\varphi_{\lambda}\in z^{\lambda}+\sum_{\beta<_r\lambda}{\cal 
C}z^{\beta}{\rm \ and \ }
 \varphi_{\lambda}\in z^{\lambda}{\cal C}[N]\leqno(6.1)$$ for each $\lambda\in 
P^+_{\Theta}$.

\begin{lemma}
   For all  $B,B'\in {\cal B}_{\theta}$, $(\chi_1,\chi_2)\in {\bf H}\times 
{\bf H}$, and $\lambda\in P^+_{\Theta}$ 
$$\varphi^{\lambda}_{B,B'}=z^{\lambda}\bar\chi_1^{-1}\bar\chi_2\cdot (z^{-\lambda}\varphi^{\lambda}_{B,B'}).$$ 
\end{lemma}

\noindent
{\bf Proof:} Let $(\chi_1,\chi_2)\in {\bf H}\times 
{\bf H}$. Note that there exists a positive 
integer $m$ such that   $P(\pi)$   is a subset of $\sum_i{\bf 
Z}(\alpha_i/m)$.
Hence for each $i$, the map  $\bar \chi_i$ can be extended  to  
a  (not necessarily unique) homomorphism, which we also call 
$\bar\chi_i$, 
  from $\sum_j{\bf 
Z}(\alpha_j/m)$ 
to  ${\cal C}^{\times}$.  
 
Consider $\lambda\in P^+(\pi)$ and fix $i\in\{1,2\}$. 
The map $\bar\chi_i$  induces  linear transformations  on $L(\lambda)$ 
and $L(\lambda)^*$,  both denoted by $\hat\chi_i$,
such that $\hat\chi_i(v)=\hat\chi_i(\wt v)v$
for each weight vector $v$ in $L(\lambda)$ or $L(\lambda)^*$.   
Suppose that $v\in L(\lambda)$ and $v'\in L(\lambda)^*$ are weight 
vectors.  Note that  $$\chi_i(u)\hat\chi_i(v)=\hat\chi_i(uv){\rm \ and \ 
}\hat\chi_i^{-1}(v') \chi_i(u)=\hat\chi_i^{-1}(v' u)
\leqno{(6.2)}$$ for all $u\in U$.  It follows that
$$c^{\lambda}_{\hat\chi_i^{-1}(v'),v}=\bar\chi_i({\rm wt}\ 
v')(c^{\lambda}_{v',v}){\rm\quad and\quad  
}c^{\lambda}_{v',\hat\chi_i(v)}=\bar\chi_i({\rm wt}\ 
v)(c^{\lambda}_{v',v}).$$
Thus 
$$\Upsilon(c^{\lambda}_{\hat\chi_1^{-1}(v'),\hat\chi_2(v)})=\bar\chi_1^{-1}({\rm wt} v')
\bar\chi_2({\rm wt}\ 
v)\Upsilon(c^{\lambda}_{v',v})=(\bar\chi_1^{-1} 
\bar\chi_2)\cdot( \Upsilon(c^{\lambda}_{v',v})).\leqno(6.3)$$

Let
$\xi$ denote a $B'$ invariant vector of $L(\lambda)$ and let $\xi^*$ 
denote a $B$ invariant vector of $L(\lambda)^*$.
Assertion  (6.2) ensures that $\hat\chi_2(\xi)$ is a nonzero $\chi_2(B')$ invariant element 
of $L(\lambda)$ and $\hat\chi_1^{-1}(\xi^*)$ is a $\chi_1(B)$ invariant vector of $L(\lambda)^*$.  
Thus
  $c_{\chi_1(B),\chi_2 
 (B')}^{\lambda}=c^{\lambda}_{\hat\chi_1^{-1}(\xi^*),\hat\chi_2(\xi)} $ up to a nonzero 
  scalar.   It follows from assertion (6.3) that $\varphi_{\chi_1(B),\chi_2 
 (B')}^{\lambda}=(\bar\chi_1^{-1} 
\bar\chi_2)\cdot(\varphi^{\lambda}_{B,B'} )$ up to a nonzero scalar.   Now $\bar\chi_1^{-1} 
\bar\chi_2(\varphi^{\lambda}_{B,B'} )$ is just a nonzero scalar 
multiple of $z^{\lambda}(\bar\chi_1^{-1} 
\bar\chi_2)\cdot(z^{-\lambda}\varphi^{\lambda}_{B,B'} )$.  The lemma follows 
from the fact that $z^{\lambda}(\bar\chi_1^{-1} 
\bar\chi_2)\cdot(z^{-\lambda}\varphi^{\lambda}_{B,B'} )$ satisfies the 
conditions of  (6.1).
$\Box$

\medskip
It follows from Lemma 6.2, that the   kernel of the action of ${\bf H}\times {\bf H}$ 
on ${\cal F}_{\theta}$ is contained in  ${\rm Ker}(\Psi)$ where 
$\Psi$ is the homomorphism of Lemma 6.1. Thus the 
action of ${\bf H}\times {\bf H}$ 
on ${\cal F}_{\theta}$
induces an action of ${\rm Hom}(N,{\cal C}^{\times})$ on ${\cal 
F}_{\theta}$. Furthermore, this action is given by 
$$g\cdot  \varphi_{\lambda}=z^{\lambda}g(z^{-\lambda}\varphi_{\lambda})
$$
for all $\lambda\in P^+_{\Theta}$,
$g\in {\rm Hom}(N,{\cal C}^{x})$, and $\{\varphi_{\lambda}|\lambda\in 
P^+_{\Theta}\}\in {\cal F}_{\theta}$.
The following result is an immediate consequence of Lemma 6.1 and 
Lemma 6.2.

\begin{theorem} The map $$(B,B')\rightarrow 
\{\varphi^{\lambda}_{B,B'}|\lambda\in P^+_{\Theta}\}$$ from ${\cal 
B}_{\theta}\times {\cal B}_{\theta}$ to ${\cal F}_{\theta}$  is a ${\bf 
H}\times {\bf H}$ equivariant map which induces a 
surjection from the set of  orbits of ${\cal B}_{\theta}\times 
{\cal B}_{\theta}$ under the action of ${\bf H}\times {\bf H}$ onto the 
set of orbits of zonal spherical families under the action of ${\rm 
Hom}(N,{\cal C}^{\times})$.
\end{theorem}

A family $\{\varphi_{\lambda}|\lambda\in P^+_{\Theta}\}\in 
{\cal F}_{\theta}$ is called $W$ invariant provided each 
$\varphi_{\lambda}\in {\cal C}[2\Sigma]^{W}$.
We next determine which   elements in ${\rm Hom}(N,{\cal 
C}^{\times})$ preserve $W$ invariance.  Recall that $Q(2\Sigma)$ is 
the root lattice associated to $2\Sigma$.

\begin{lemma}
Let $g\in {\rm Hom}(N,{\cal C}^{\times})$ and let 
$\{\varphi_{\lambda}|\lambda\in P^+_{\Theta}\}$ be a $W$ 
invariant zonal spherical family.   Then $g\cdot 
\varphi_{\lambda}$
is $W$ invariant for all $\lambda\in P^+_{\Theta}$  
if and only if $g$ 
acts trivially on 
    $Q(2\Sigma)$.  
     \end{lemma}
    
\noindent
{\bf Proof:} Suppose that $\Sigma$ is of type $BC_r$.   Then
the weight lattice of $\Sigma$ is equal to the weight lattice of the 
subroot system 
$\Sigma_1$ in $\Sigma$ of type $B_r$. In particular, if the lemma holds when $\Sigma$ 
is of type $B_r$, then it follows for $\Sigma$ of type $BC_r$.  Hence 
we can reduce to the case where $\Sigma$ is a reduced root system.    

Consider $g\in {\rm Hom}(N, {\cal C}^{\times})$. Note that $N$ is a free 
${\bf Z}$ module.  Furthermore, there exists a positive integer $m$ 
such that $P(\Sigma)$ is a sublattice of $\{\beta/m|\beta\in N\}$.
 Thus we may extend
$g$ to a (not necessarily unique) element of ${\rm 
Hom}(P(\Sigma),{\cal C}^{\times})$.   Hence, it is sufficient to prove the 
assertions of the theorem with ${\rm Hom}(N, {\cal C}^{\times})$ replaced by 
${\rm 
Hom}(P(\Sigma),{\cal C}^{\times})$.

Recall that $\{\varphi_{\lambda}|\lambda\in P^+_{\Theta}\}$ is a 
basis for ${\cal C}[2\Sigma]^{W}$. Thus, the first 
asssertion of the lemma is equivalent to  
  $g$ restricts to an action on  ${\cal C}[2\Sigma]^{W}$  if and 
  only if $g$ acts trivially on $Q(2\Sigma)$.

Assume first that   $g$ sends ${\cal C}[2\Sigma]^{W}$ to 
itself. Given $X=\sum_{\beta}a_{\beta}z^{\beta}$, set ${\rm 
Supp}(X)=\{z^{\beta}|a_{\beta}\neq 0\}$. 
Let $\delta$ denote the half sum of the positive restricted roots in $\Sigma$.
By [H, Corollary 10.2], 
$2(\delta,\tilde\alpha_i)=(\tilde\alpha_i,\tilde\alpha_i)$ and 
$\tilde s_i(\delta)=\delta-\tilde\alpha_i$ for each 
$\tilde\alpha_i\in\Sigma$.  It follows that $2\delta\in P(2\Sigma)$. 
  By say [H, 
Theorem 10.3(e)], $w\delta=\delta$ if and only if $w=1$ and thus 
$w\delta=w'\delta$ if and only if $w=w'$. Thus for each $i$ such 
that $\tilde\alpha_i\in \Sigma$, 
 we can write
$$m_{\delta}=z^{2\delta}+z^{2\delta-2\tilde\alpha_i}+\sum_{\beta\neq 
2\delta,\beta\neq 
2\delta-2\tilde\alpha_i}a_{\beta}z^{\beta}\leqno{(6.4)}$$  for 
some scalars $a_{\beta}$.

   By  (6.4), we have
 $$g\cdot m_{\delta}=g({2\delta})z^{\delta}+ 
 g({2\delta-2\tilde\alpha_i})z^{\delta-\tilde\alpha_i}+\sum_{\beta\neq 
\delta,\beta\neq 
\delta-\tilde\alpha_i}g({\beta})a_{\beta}z^{\beta}.\leqno{(6.5)}$$
Note that ${\rm Supp}(g\cdot m_{\delta})={\rm Supp}(m_{\delta})$
and $g\cdot m_{\delta}\in 
{\cal C}[2\Sigma]^{W}$.  Recall that the 
$m_{\beta},$ for $\beta\in P^+_{\Theta}$, form a basis for ${\cal C}[2\Sigma]^{W}$.  
Furthermore, by [H, Section 13.2, Lemma A],
the sets 
${\rm Supp}(m_{\beta})$ are pairwise disjoint.   Hence  (6.5) implies 
that $g\cdot m_{\delta}$ is 
equal to $g({2\delta})m_{\delta}$.  It further follows from (6.5) that
 $g({2\delta})=g({2\delta-2\tilde\alpha_i})=g({2\delta})g({2\tilde\alpha_i})^{-1}$.
Therefore $g({2\tilde\alpha_i})=1$ for all $i$.  Hence $g$ acts 
trivially on $ Q(2\Sigma)$. 

Now assume that $g\in {\rm Hom}(P(\Sigma),{\cal C}^{\times})$ such that $g({2\tilde\alpha_i})=1$ for 
all $\alpha_i\in \pi-\pi_{\Theta}$.
 Note that the set $2\Sigma$ is a root system with respect to the 
Cartan inner product (of the same type as 
$\Sigma$) with set of positive 
simple roots equal to $\{2\tilde\alpha_i|\alpha_i\in\pi-\pi_{\Theta}\}$.
Since we are assuming that $\Sigma$ is  reduced, we can 
find a semisimple Lie algebra $L$ with root system equal to 
$2\Sigma$.  It is well known  that the characters of the finite dimensional simple 
modules associated to $L$ are parametrized by the dominant integral 
weights $P^+_{\Theta}$  and  form a basis for ${\cal 
C}[2\Sigma]^{W}$. Let 
$X_{\lambda}$ be the character associated to the finite dimensional 
simple $L$ module $V(\lambda)$ of highest weight $\lambda$ for $\lambda\in 
P^+_{\Theta}$.  Now the  root vectors of $L$ have weight $\beta$, 
for $\beta\in 2\Sigma$. Set 
 $\lambda+Q(2\Sigma)=\{\lambda+\beta|\beta\in Q(2\Sigma)\}$.
It follows that the weights of $V(\lambda)$ are contained in 
the set $\lambda+Q(2\Sigma)$.
  In particular, ${\rm 
Supp}(X_{\lambda})\subset \lambda+Q(2\Sigma)$.    Since $g$ 
acts as the identity on $Q(2\Sigma)$, it follows that 
$g\cdot z^{\lambda-\beta}=g({\lambda})z^{\lambda-\beta}$ for all 
$\beta\in Q(2\Sigma)$.  Therefore $g\cdot X_{\lambda}$ is equal to 
$g({\lambda})X_{\lambda}$ for each $\lambda\in P^+_{\Theta}$. Thus
 $g$ sends ${\cal C}[2\Sigma]^{W}$ to itself.
$\Box$ 

\medskip
Set ${\cal G}_{\theta}$ equal to the subgroup of      
 ${\rm Hom}(N, {\cal C}^{\times})$ which acts trivially on $Q(2\Sigma)$.  In 
 particular, ${\cal 
G}_{\theta}$ is 
isomorphic to ${\rm Hom}(N/Q(2\Sigma),{\cal C}^{\times})$.   Note that 
$N/Q(2\Sigma)$ is isomorphic to the direct product of 
$|{\cal S}|$ copies of ${\bf Z}_2$.  Thus ${\cal G}_{\theta}\cong
{\rm Hom}(N/Q(2\Sigma),\{1,-1\})\cong {\bf Z}_2^{|{\cal S}|}.$

 Recall that
 $\{B_{\theta,{\bf s},{\bf d}}|{\bf s}\in {\bf S}({\cal C})$ 
and ${\bf d}\in {\bf D}({\cal C})\}$ is a complete set of distinct 
representatives for the orbits of  ${\cal B}_{\theta}$ under the 
action of ${\bf H}$.  Set ${\cal O}(B_{\theta,{\bf s},{\bf d}},
B_{\theta,{\bf s}',{\bf d}'})$ equal to the orbit containing $ (B_{\theta,{\bf s},{\bf d}},
B_{\theta,{\bf s}',{\bf d}'})$  under the action of ${\bf 
H}\times {\bf H}$.  Note that  $\{{\cal O}(B_{\theta,{\bf s},{\bf d}},
B_{\theta,{\bf s}',{\bf d}}')|{\bf s},{\bf s}'\in {\bf S}({\cal C})$ 
and ${\bf d},{\bf d}'\in {\bf D}({\cal C})\}$ is a complete list  of 
the  
distinct  orbits of ${\cal B}_{\theta}\times 
{\cal B}_{\theta}$ under this action.
 Given ${\bf s},{\bf s}'\in {\bf S}({\cal C})$ and ${\bf d},{\bf d}'\in 
 {\bf D}({\cal C})$, 
set $$\varphi^{\lambda}_{{\bf s},{\bf s}',{\bf d}, {\bf d}'}= \varphi^{\lambda}_{B,B'}$$ where
$B=B_{\theta,{\bf s}, {\bf d}}$,  $B'=\chi_{\bf c}(B_{\theta,{\bf 
s}',{\bf d}'})$, and  ${\bf c}\in {\bf D}({\cal C})$ 
is chosen to satisfy ${\bf c}^2={\bf 
d}^{-1}{\bf d}'$.

\begin{theorem}  The set 
$\{\varphi_{\lambda}|\lambda\in P^+_{\Theta}\}$  is a $W$ 
invariant zonal spherical family associated to an 
element of ${\cal O}(B_{\theta,{\bf s},{\bf d}},
B_{\theta,{\bf s}',{\bf d}'})$
 if and only if   there exists $g\in {\cal G}_{\theta}$ 
such that $$\varphi_{\lambda}=g\cdot \varphi^{\lambda}_{{\bf s},{\bf 
s}',{\bf d}, {\bf d}'}\leqno{(6.6)}$$ for all $\lambda\in 
P^+_{\Theta}$. 
In particular, if ${\cal S}$ is empty, then  $g$ in 
(6.6) is just the identity map.
\end{theorem}

\noindent
{\bf Proof:}   
  Let $B= B_{\theta,{\bf s},{\bf d}}$ and  $B'=\chi_{\bf c}(B_{\theta,{\bf 
s}',{\bf d}'})$, where  ${\bf c}$ is the element of ${\bf D}$ which satisfies ${\bf c}^2={\bf 
d}^{-1}{\bf d}'$ chosen  above.  By Corollary 5.4, $\{\varphi^{\lambda}_{{\bf s},{\bf 
s}',{\bf d}, {\bf d}'}|\lambda\in P^+_{\Theta}\}$ is a 
a  $W$ invariant zonal spherical family associated to $(B, 
B')$.  Let 
$\{\varphi_{\lambda}|\lambda\in P^+_{\lambda}\}$ be another 
zonal spherical family associated to an element in ${\cal O}_1$.  
By Theorem 6.3, there exists $g\in {\rm Hom}(N,{\cal C}^{\times})$ such that 
$$\{\varphi_{\lambda}|\lambda\in P^+_{\Theta}\}=g\cdot \{\varphi ^{\lambda}_{{\bf s},{\bf 
s}',{\bf d}, {\bf d}'}|\lambda\in P^+_{\Theta}\}.$$ The result now 
follows from Lemma 6.4. 
 $\Box$
  
  \medskip
  
  Consider the case when $\gg,\theta$ is an irreducible pair.  There 
  are three possibilities:
  \begin{enumerate} 
  \item[(i)] ${\cal S}$ and ${\cal D}$ are both empty.
  \item[(ii)] ${\cal S}$ is empty and ${\cal D}$ has exactly one element.
  \item[(iii)] ${\cal S}$ has exactly one element and ${\cal D}$ is empty.
  \end{enumerate}
 
Note that when ${\cal S}$ is empty, ${\cal G}_{\theta}$ is the trivial 
group.  In the first case, ${\cal B}_{\theta}$ 
is a single orbit under the action of ${\bf H}$.  Thus 
Theorem 6.5 implies that we can associate a unique $W$ invariant zonal 
spherical family   to ${\cal B}_{\theta}$ when both ${\cal S}$ and 
${\cal D}$ are empty. Suppose ${\cal D}=\{\alpha_i\}$  as in 
case (ii).  Then the set ${\bf D}({\cal C})$ is a  set 
of $n$-tuples which vary only in the entry $d_i$.  Thus, 
there is a unique two-parameter set of $W$ invariant zonal 
spherical families associated   to ${\cal B}_{\theta}$. If ${\cal S}$ 
and ${\cal D}$ satisfy condition (iii) then
${\cal G}_{\theta}$ is a cyclic group of order $2$.  In this case, we can 
associate exactly
two two-parameter sets of $W$ invariant zonal 
spherical families    to ${\cal B}_{\theta}$.

It should be noted that the existence of   ``two-parameter'' families of 
$W$ zonal 
spherical functions associated to a quantum symmetric pair has appeared
 in the literature in a few special cases.  Indeed, 
  two-parameter families of zonal sperical fuctions are studied in 
  [K1] when $\gg={\bf sl}\ 2$ and 
 $\gg,\theta$ is of type ${\bf AI}$ and in [DK] when $\gg,\theta$ is 
 of type ${\bf AIV}$. (Types ${\bf AI}$ and ${\bf AIV}$ are explained 
 in the next section.)

\section{ Irreducible Symmetric Pairs}
 By  
[A, 2.5 and 5.1],   the pair 
$\gg,\theta$ is   irreducible if and only if 
\begin{enumerate}
\item[(7.1)] $\gg$ is simple \ \ \    or
\item[(7.2)] $\gg= \gg_1\oplus \gg_2$ where both $\gg_1$ and $\gg_2$ 
are simple. Here  $\gg_1$ is generated by 
 $\{e_i,f_i,h_i|1\leq i\leq m\}$ and $\gg_1$ is generated by 
 $\{e_{i+m},f_{i+m},h_{i+m}|1\leq i\leq m\}$.  Moreover,
 $$e_i\rightarrow e_{i+m},f_i\rightarrow f_{i+m},h_i\rightarrow 
 h_{i+m}$$ defines an isomorphism from $\gg_1$ to $\gg_2$, and  
 $\theta$ is defined by 
 $$\theta(e_i)=f_{i+m},\quad \theta(f_i)=e_{i+m},\quad 
 \theta(h_i)=h_{i+m}.$$
\end{enumerate}

Assume for the moment that $\gg$ and $\theta$ are as described in (7.2). 
In this case, both ${\cal D}$ and ${\cal 
S}$ are empty sets. The group $T_{\Theta}$ is generated by 
$t_it_{i+m}^{-1}$, $1\leq i\leq m$ and the algebra 
$B_{\theta,{\bf s},{\bf d}}$ is generated by 
$T_{\Theta}$ together with  the 
elements   
$B_i=y_it_i+t_{i+m}^{-1}x_{i+m}t_i$ and
 $B_{i+m}=y_{i+m}t_{i+m}+ t_{i}^{-1}x_{i}t_i$ for $1\leq i\leq m$.

 Now consider the case when $\gg$ is simple.   
  We use 
the classification in 
[A] (see also [He, Chapter X, Section F]) to give 
a complete list of the subalgebras 
$B_{\theta,{\bf s},{\bf d}}$ associated to  
 irreducible pairs 
$\gg,\theta$ for $\gg$ simple. In each case, we describe 
$\pi_{\Theta}$, ${\cal S}$, ${\cal D}$, the permutation $p$ when $p$ 
is not the identity, and  $B_i, 
\alpha_i\notin\pi_{\Theta}$.  This information 
is enough to determine  all the generators of $B_{\theta,{\bf s},{\bf 
d}}$ since ${\cal M}$ is just the quantized enveloping 
algebra associated to the root system
$\pi_{\Theta}$ considered as a subset of $\pi$ and $T_{\Theta}$ is 
generated by the sets $\{t_i|  \alpha_i\in\pi_{\Theta}\}$ and 
$\{t_it_{p(i)}^{-1}|
\alpha_i\notin\pi_{\Theta}\}$. 

We use below the numbering of the vertices of the Dynkin diagram 
given in [H, Section 11.4].   It should be noted that when $\gg$ is of 
type $D_n$, the roots $\alpha_n$ and $\alpha_{n-1}$ can be reordered 
with $n$ and $n-1$ interchanged.

\medskip
\noindent
{\bf Type AI}
 $\gg$ is of type $A_n$, $\pi_{\Theta}$ is the empty set, 
${\cal S}$ and ${\cal D}$ are both empty,  $B_i=y_it_i+q_i^{-2}x_i$ for each 
$1\leq i\leq n$.

\medskip
\noindent
{\bf Type AII}
 $\gg$ is of type $A_n$, where $n=2m+1$ is odd and $n\geq 3$, 
$\pi_{\Theta}= \{\alpha_{2j+1}| 0\leq j\leq m\}$, 
  ${\cal S}$ and ${\cal D}$ are both empty, 
$B_i=y_it_i+[(\adr x_{i-1} x_{i+1})t_i^{-1}x_i]t_i$,
for $i=2j, 1\leq j\leq m$. 

\medskip
\noindent
{\bf Type AIII}
{\bf{\it Case 1}}:
 $\gg$ is of type $A_n$, $r$ is an integer such that $2\leq r\leq n/2$ 
$\pi_{\Theta}=\{\alpha_j|\ r+1\leq j\leq n-r\}$,
${\cal S}$ is empty,  ${\cal D}=\{\alpha_{r}\}$,  
  $p$   sends $i$ to $n-i+1$ for each $1\leq i\leq n$,
$$B_i=y_it_i+q_i^{-2}x_{p(i)}t_{p(i)}^{-1}t_i {\rm \ \ for\ \ }
 1\leq i\leq r-1{\rm \ \ and\ \ }n-r+2\leq i\leq n,$$ 
$$B_{r}=y_{r}t_{r}+d_r[(\adr x_{r+1})(\adr x_{r+2})\cdots(\adr 
x_{n-r})t_{p(r)}^{-1}x_{p(r)}]t_{r},$$ and 
$$B_{p(r)}=y_{p(r)}t_{p(r)}+(-1)^{n-2r}[(\adr x_{n-r})(\adr 
x_{n-r-1})\cdots(\adr 
x_{r+1})t_{r}^{-1}x_{r}]t_{p(r)}.$$

\noindent{\bf{\it Case 2}}:  $\gg$ is of type $A_n$ where $n=2m+1$, 
$\pi_{\Theta}$ is empty,  ${\cal D}$ is empty, ${\cal S}=\{\alpha_{m+1}\},$
 $p$   sends $i$ to $n-i+1$ for $1\leq i\leq n$,
$B_i=y_it_i+q_i^{-2}x_{p(i)}t_{p(i)}^{-1}t_i$
for $1\leq i\leq m$ and $m+2\leq i\leq n$, and 
$B_{m+1}=y_{m+1}t_{m+1}+q_{m+1}^{-2}x_{m+1}+s_{m+1}t_{m+1}.$

\medskip
\noindent
{\bf Type AIV} 
Same as type AIII, case 1, with $r=1$.

\medskip
\noindent
{\bf Type BI}  
 $\gg$ is of type $B_n$, $r$ is an integer such that $2\leq r\leq n$,
$\pi_{\Theta}=\{\alpha_i|r+1\leq i\leq n\}$, 
${\cal S}$ and ${\cal D}$ 
are empty,
$B_i=y_it_i+q_i^{-2}x_i$ for $1\leq i\leq r-1$
and $${B_r=y_rt_r
+[(\adr x_{r+1}\cdots x_{n-1} x_n^{(2)} x_{n-1}\cdots  x_{r+1})t_rx_r]t_r^{-1}.}$$

\medskip
\noindent
{\bf Type BII}  
Same as type BI, Case 1, only $r=1$. 

\medskip
\noindent
{\bf Type CI}
 $\gg$ is of type $C_n$, $\pi_{\Theta}$ is empty, 
${\cal D}$ is empty,  ${\cal S}=\{\alpha_{n}\}$, $B_i=y_it_i+q_i^{-2}x_i$ for $1\leq 
i\leq n-1$, and 
$B_{n}=y_{n}t_{n}+q_n^{-2}x_{n}+s_nt_{n}.$

\medskip
\noindent
{\bf Type CII}
{\bf{\it Case 1}}:  $\gg$ is of type $C_n$, $r$ is an even integer such 
that $1\leq r\leq (n-1)$,
$\pi_{\Theta}=\{\alpha_{2j-1}|1\leq j\leq r/2\}\cup\{\alpha_j|r+1\leq 
j\leq n\}$,  
 ${\cal S}$ and ${\cal D}$ are both empty,
$$B_i=y_it_i+[(\adr x_{i-1} x_{i+1})t_i^{-1}x_i]t_i$$
for $i=2j, 1\leq j\leq (r-2)/2$
and $$B_r=y_rt_r+[(\adr x_{r-1}x_{r+1}\cdots 
x_{n-1}x_nx_{n-1}\cdots x_{r+1})t_r^{-1}x_r]t_r.$$

\medskip
\noindent
{\bf{\it Case 2}}:
 $\gg$ is of type $C_n$ where $n$ is even,
$\pi_{\Theta}=\{\alpha_{2j-1}|1\leq j\leq n/2\}$,
 ${\cal S}$ and ${\cal D}$ are both empty,
$B_i=y_it_i+(\adr x_{i-1} x_{i+1})t_i^{-1}x_i$
for $i=2j, 1\leq j\leq (n-2)/2$
and $B_n=y_nt_n+[(\adr
x_{n-1}^{(2)}) t_n^{-1}x_n]t_n.$

\medskip
\noindent
{\bf Type DI}
{\bf{\it Case 1}}:  
  $\gg$ is of type $D_n$,  $r$ is an integer such that $2\leq r\leq 
  n-2$, 
$\pi_{\Theta}=\{\alpha_i|r+1\leq i\leq n\}$, 
${\cal S}$ and ${\cal D}$ 
are both empty,
$B_i=y_it_i+q_i^{-2}x_i$ for $1\leq i\leq r-1$,
and $${B_r=y_rt_r
+[(\adr x_{r+1}\cdots x_{n-2}
 x_{n-1}
x_nx_{n-2}\cdots  x_{r+1})t_r^{-1}x_r]t_r.}$$

\medskip
\noindent
{\bf{\it Case 2}}: We assume $n\geq 4$ (the case when $n=3$ is the same as 
type AI.) $\gg$ is of type $D_n$,  
$\pi_{\Theta}$ is empty, 
${\cal S}$ and ${\cal D}$ 
are both empty,  
$p(i)=i$ for $1\leq i\leq n-2$, $p(n-1)=n$, and $p(n)=n-1$  and 
$B_i=y_it_i+t_{p(i)}^{-1}x_{p(i)}t_i$ for $1\leq i\leq n.$

\medskip
\noindent
{\bf{\it Case 3}}: $\gg$ is of type $D_n$,  
$\pi_{\Theta}$ is empty, 
${\cal S}$ and ${\cal D}$ 
are both empty, and
$B_i=y_it_i+q_i^{-2}x_i$ for $1\leq i\leq n$.

\medskip
\noindent
{\bf Type DII} 
This is the same as DI, Case 1, with $r=1$.

\medskip
\noindent
{\bf Type DIII}
{\bf{\it Case 1}}:
 $\gg$ is of type $D_n$ where $n$ is even, 
$\pi_{\Theta}=\{\alpha_{2i-1}|1\leq 
i\leq 
n/2\}$, ${\cal D}$ is empty,  ${\cal S}=\{\alpha_{n}\}$,
$B_i=y_it_i+(\adr x_{i-1} x_{i+1})t_i^{-1}x_i$
for $i=2j, 1\leq j\leq (n-2)/2$,
and 
$B_n=y_nt_n+q_n^{-2}x_n+s_nt_n.$
  
\medskip
\noindent
{\bf{\it Case 2}}:  $\gg$ is of type $D_n$, $n$ is odd, 
$\pi_{\Theta}=\{\alpha_{2i-1}|1\leq 
i\leq 
(n-1)/2\}$, ${\cal S}$ is empty, ${\cal D}=\{\alpha_{n-1}\}$,
$B_i=y_it_i+(\adr x_{i-1} x_{i+1})t_i^{-1}x_i$
for $i=2j, 1\leq j\leq (n-3)/2$, $p(j)=j$ for $1\leq j\leq n-2$,
$p(n)=n-1$, and $p(n-1)=p(n)$,
$B_{n-1}=y_{n-1}t_{n-1}+d_{n-1}[(\adr x_{n-2})t_{n}^{-1}x_{n}]t_n$
and 
$B_n=y_nt_n-[(\adr x_{n-2})t_{n-1}^{-1}x_{n-1}]t_n.$

\medskip
\noindent
{\bf Type E1, EV, EVIII}
 $\gg$ is of type E6,E7, E8 respectively, $\pi_{\Theta}$ 
is empty,   ${\cal S}$ and 
${\cal D}$ are both empty, and  $B_i=y_it_i+q_i^{-2}x_i$,
for all $\alpha_i\in \pi$.

\medskip
\noindent
{\bf Type EII}
 $\gg$ is of type E6, $\pi_{\Theta}$ is empty, 
 both ${\cal S}$ and ${\cal 
D}$ are empty,  
$p(1)=6,p(3)=5,p(4)=4,p(2)=2,p(5)=3,$ and $p(6)=1$,
and
$B_i=y_it_i+q_i^{-2}x_{p(i)}t_{p(i)}^{-1}t_i$
for $1\leq i\leq 6$.

\medskip
\noindent
{\bf Type EIII} 
 $\gg$ is of type E6, 
$\pi_{\Theta}=\{\alpha_3,\alpha_4,\alpha_5\}$,
${\cal S}$ is empty, ${\cal D}=\{\alpha_1\}$,  
$p(1)=6,p(3)=5,p(4)=4,p(2)=2,p(5)=3,p(6)=1$,
$$B_1= y_1t_1+d_1[(\adr x_3)(\adr x_4)(\adr x_5)t_6^{-1}x_6]t_1,$$
$$B_6= y_6t_6+[(\adr x_5)(\adr x_4)(\adr x_3)t_1^{-1}x_1]t_6,$$
and 
$$B_2=y_2t_2+[(\adr x_4x_3x_5x_4)t_2^{-1}x_2]t_2.$$

\medskip
\noindent
{\bf Type EIV}
 $\gg$ is of type E6, 
$\pi_{\Theta}=\{\alpha_3,\alpha_4,\alpha_5,\alpha_2\}$,
${\cal S}$ and ${\cal D}$ are both  empty,
$$B_1=y_1t_1+[(\adr x_3x_4x_5x_2x_4x_3)t_1^{-1}x_1]t_1$$ and 
$$B_6=y_6t_6+[(\adr x_5x_4x_3x_2x_4x_5)t_6^{-1}x_6]t_6.$$

\medskip
\noindent
{\bf Type EVI}  $\gg$ is of type E7, 
$\pi_{\Theta}=\{\alpha_7,\alpha_5,\alpha_2\}$, both ${\cal S}$ 
and ${\cal D}$ are empty, $B_6=y_6t_6+[(\adr x_7x_5)t_6^{-1}x_6]t_6,$
$B_4=y_4t_4+[(\adr x_2x_5)t_4^{-1}x_4]t_4,$
and $B_i=y_it_i+q_i^{-2}x_i$ for $i=1,3$.

\medskip
\noindent
{\bf Type EVII} 
 $\gg$ is of type E7,
$\pi_{\Theta}=\{\alpha_2,\alpha_3,\alpha_4,\alpha_5\}$, ${\cal D}$ is 
empty, ${\cal S}=\{\alpha_7\}$,  $$B_1=y_1t_1+[(\adr x_3x_4x_2x_5x_4x_3)t_1^{-1}x_1]t_1,$$
$$B_6=y_6t_6+[(\adr x_5x_4x_2x_3x_4x_5)t_6^{-1}x_6]t_6,$$
and $B_7=y_7t_7+q_7^{-2}x_7+s_7t_7.$

\medskip
\noindent
{\bf Type EIX}
$\gg$ is of type E8,
$\pi_{\Theta}=\{\alpha_2,\alpha_3,\alpha_4,\alpha_5\}$, ${\cal D}$ and
 ${\cal S}$ are both empty, $$B_1=y_1t_1+[(\adr x_3x_4x_2x_5x_4x_3)t_1^{-1}x_1]t_1,$$
$$B_6=y_6t_6+[(\adr x_5x_4x_2x_3x_4x_5)t_6^{-1}x_6]t_6,$$
and $B_i=y_it_i+q_i^{-2}x_i$ for $i=7,8.$

\medskip
\noindent
{\bf Type FI} 
 $\gg$ is of type F4, $\pi_{\Theta}$ is empty,
${\cal S}$ and ${\cal D}$ are empty,
$B_i=y_it_i+q_i^{-2}x_i$ for $i=1,2,3,4$.

\medskip
\noindent
{\bf Type FII}
 $\gg$ is of type F4, 
$\pi_{\Theta}=\{\alpha_1,\alpha_2,\alpha_3\}$,
both ${\cal D}$ and ${\cal S}$ are  empty,  and 
$$B_4=y_4t_4+[(\adr x_3x_2x_1x_3x_2x_3)t_4^{-1}x_4]t_4$$

\medskip
\noindent
{\bf Type G}
$\gg$ is of type G2, $\pi_{\Theta}$ is empty,  ${\cal 
S}$ and ${\cal D}$ are both empty, and 
$B_i=y_it_i+q_i^{-2}x_i$ for $i=1,2$.

\medskip

Set $B_i=y_it_i$ for $\alpha_i\in \pi_{\Theta}$. 
Recall that ${\cal M}^+={\cal M}\cap U^+$.
 Given $\alpha_i\notin\pi_{\Theta}$, let $Z_i$ be the element in ${\cal M}^+$ such 
 that $$\tilde\theta(y_i)=(\adr Z_i)t_{p(i)}^{-1}x_{p(i)}.$$   
 Theorem 7.4 of [L4] 
(see also [L4, Variations 1 and 2]) gives most of the relations satisfied by the generators 
 of $B_{\theta,{\bf s}, {\bf d}}$ 
 and a procedure for computing the remaining ones. (This proof of 
 this result is 
 stated in [L4] 
 for the $B_{\theta,{\bf s},{\bf d}}$ where ${\bf 
 s}\in {\bf S}({{\bf C}[q]}_{(q-1)})$ and ${\bf d}\in {\bf D}({{\bf 
 C}[q]}_{(q-1)})$, but  the arguments
 work for   all 
 $B_{\theta,{\bf s}, {\bf d}}$,
 for ${\bf s}\in {\bf S}({\cal C})$ and ${\bf d}\in {\bf D}({\cal C})$.)
 Using the above complete list of the left coideal subalgebras 
 $B_{\theta,{\bf s}, {\bf d}}$,
 we make these relations  precise. 
  In order to make the notation easier in the following theorem, 
we define a modified version $p'$ of the permutation  $p$    such that 
 $p'=p$ on the set $\{i|\alpha_i\notin\pi_{\Theta}\}$ and $p'(j)=0$ 
for all $j$ in the set $\{j|\alpha_j\in \pi_{\Theta}\}$.

 \begin{theorem} The left coideal subalgebra $B_{\theta,{\bf s},{\bf 
 d}} $
 of $U$ is the algebra 
 generated  over 
${\cal M}^+T_{\Theta}$ by the elements   $ B_i$, $1\leq i\leq 
n$ 
subject to the following relations:
 \begin{enumerate}
 \item[(i)] $\tau(\lambda) 
  B_i\tau(-\lambda)=q^{-(\lambda,\alpha_i)}   B_i$ for all
$\tau(\lambda)\in T_{\Theta}$ and $1\leq i\leq n$. 
\item[(ii)] $ t_j^{-1} x_j   B_i-   B_i
t_j^{-1} x_j=\delta_{ij}(t_j-t_j^{-1})/(q_j-q_j^{-1})$ for all 
$\alpha_j\in \pi_{\Theta}$ and $1\leq i\leq n.$  
 \item[(iii)] 
 If $a_{ij}=0$, then  $$\eqalign{ &B_i  B_j-  B_j  B_i = \cr 
  &\delta_{p'(i),j}(q_i-q_i^{-1})^{-1}(-d_i[(\adr 
  Z_i)t_{p(i)}^{-2}]t_it_{p(i)}-d_{p(i)}[(\adr 
  Z_{p(i)})t_i^{-2}]t_it_{p(i)}).\cr}$$
\item[(iv)] If $a_{ij}=-1$, then 
$$\eqalign{& B_i^2  B_j-(q_i+q_i^{-1})  B_i 
B_j B_i+  B_j  B_i^2=\cr
& \delta_{i,p'(i)}q_i^{-1}
[(\adr 
Z_i)t_i^{-2}]t_i^2 B_j 
 -\delta_{i,p'(j)}(q_i+q_i^{-1}) B_i(d_iq_i^{-1} 
t_j^{-1}t_i+d_jq_i^{2} t_i^{-1}t_j).\cr}$$
\item[(v)] If $a_{ij}=-2$, then
$$\eqalign{& B_i^3  B_j+(q_i^2+1+q_i^{-2})(- 
B_i^2  B_j B_i+  B_i  B_j  B_i^2)- 
B_j B_i^3=\cr
&
\delta_{i,p'(i)}q_i^{-1}(q_i+q_i^{-1})^2(  B_i  B_j- 
B_j  B_i).\cr}$$
\item[(vi)] If $a_{ij}=-3$, then 
$$\eqalign{& B_i^4 B_j+B_j
B_i^4-(q_i^3+q_i+q_i^{-1}+q_i^{-3}) B_i^3B_j
B_i + B_i B_j B_i^3)\cr
&+(q_i^4+q_i^2+2+q_i^{-2}+q_i^{-4}) B_i^2 B_j
B_i^2 \cr& 
= (q_i^{-5}+2q_i^{-3}+4q_i^{-1}+2q_i+q_i^{3})( B_i^2
 B_j+ B_j B_i^2)\cr 
&-(q_i^4+4q_i^2+5+5q_i^{-2}+4q_i^{-4}+q_i^{-6}) B_i
B_j B_i
-(q_i^2+1+q_i^{-2}) B_j. \cr}$$
  \end{enumerate}
 \end{theorem}

The relations in Theorem 7.1 can be checked directly in a case-by-case fashion 
 using the list of algebras $B_{\theta,{\bf s},{\bf d}}$ above.  
Indeed, this is the approach used in [L1] which handles most of the 
cases when
 $\Theta(\alpha_i)=-\alpha_{p(i)}$ 
 for $1\leq i\leq n$. Similar 
 computations can be made in the other cases. However, in 
 general, these computations are rather brutal.  In this paper,
  we  use the alternate, quicker method derived from the proof of 
  [L2, Theorem 7.4]. 
   First, we 
 review some notation from [L4] and prove a small technical lemma.
 
 Let $G^+$ be the subalgebra of $U$ generated by 
 $x_it_i^{-1},1\leq i\leq n$. As in [L4, (4.2)]. we have the following two direct sum 
 decompositions  $$U=\sum_{\lambda,\mu}U^-_{-\lambda}G^+_{\mu}U^o{\rm \ 
 and \ }U^o=\sum_{\lambda}{\cal
 C}\tau(\lambda).$$  Set $\pi_{\lambda,\mu}$ equal to the 
 projection of $U$ onto $U^-_{-\lambda}G^+_{\mu}U^o$ using the above 
 decomposition of $U$ and
 set $P_{\gamma}$ equal to the projection of $U^o$ onto ${\cal
 C}\tau(\gamma)$ using the decomposition of $U^o$.

 Let $Z$ be an element of ${\cal M}^+$. 
 The form of the right adjoint action (1.2) implies that   
 $[(\adr Z)t_{p(i)}^{-2}]$ is an element of 
 ${\cal M}t_{p(i)}^{-2}$.   Since $t_{p(i)}^{-1}t_i$ is in 
 $T_{\Theta}$, it follows that $[(\adr Z)t_{p(i)}^{-2}]t_{p(i)}t_i$ 
 is an element of $B_{\theta,{\bf s},{\bf d}}$.

\begin{lemma} Let  $Z\in {\cal M}^+$ be a vector of weight $\gamma$.
 For all 
$\alpha_j\in\pi_{\Theta}$ and each $\mu\in Q^+(\pi)$, $$\pi_{\alpha_j,\mu}\circ
{\it\Delta}((\adr Z)t_j^{-1}x_j)=\delta_{\mu,\gamma}t_j^{-1}x_j\otimes 
[(\adr Z)t_j^{-2}]t_j.$$
\end{lemma}

\noindent
{\bf Proof:}     We 
show how the argument works when $\gamma$ is a simple root. The 
general case follows similarly using induction on the weight of $Z$.  
Thus, assume that $Z=x_i$ for some 
$\alpha_i\in\pi_{\Theta}$.   Recall ([Jo, 3.2.9]) that 
$${\it\Delta}(x_k)=x_k\otimes 1+t_k\otimes x_k$$ for each $1\leq 
k\leq n$. Hence (1.2) implies that  $$\eqalign{&{\it\Delta}((\adr 
x_i)(t_j^{-1}x_j))\cr&=
-(x_it_i^{-1}\otimes t_i^{-1}+1\otimes x_it_i^{-1})(t_j^{-1}x_j\otimes 
t_j^{-1}+1\otimes t_j^{-1}x_j)\cr&+(t_i^{-1}\otimes 
t_i^{-1})(t_j^{-1}x_j\otimes t_j^{-1}+1\otimes 
t_j^{-1}x_j)(x_i\otimes 1+t_i\otimes x_i).\cr}\leqno{(7.1)}$$
Inspection of (7.1) shows that  $\pi_{\alpha_j,\mu}\circ
{\it\Delta}((\adr x_i)t_j^{-1}x_j)=0$ unless $\mu=\alpha_i$. Thus, 
a straightforward computation using (7.1) yields  $$\eqalign{\pi_{\alpha_j,\alpha_i}\circ
{\it\Delta}((\adr x_i)t_j^{-1}x_j)&=t_j^{-1}x_j\otimes 
(-x_it_i^{-1}t_j^{-1}
+q^{-(\alpha_i,\alpha_j)}t_i^{-1}t_j^{-1}x_i)\cr&=
t_j^{-1}x_j\otimes [(\adr x_i)t_j^{-2}]t_j.\Box\cr}$$

 \noindent
 {\bf Proof of Theorem 7.1:}  Relations (i) and (ii) are just (i) and 
 (ii) of [L4, 
 Theorem 7.4]. We check relations (iii) through (vi).

Given $B_i$ and $B_j$, set $$Y=Y(B_i,B_j)=
\sum_{m=0}^{1-a_{ij}}(-1)^m\left[\eqalign{1&-a_{ij}\cr 
 &m\cr}\right]_{q_i} B_i^{1-a_{ij}-m} B_j B_i^{m}.$$
 Set $\lambda= (1-a_{ij})\alpha_i+\alpha_j$.
 By the proof of [L4, Theorem 7.4] (see also [L4, Variations  1 and 2] ), $$0=((P_{\lambda}\circ\pi_{0,0})\otimes 
 Id){\it\Delta}(Y).\leqno(7.2)$$  Let $Y'\in B$ such that 
 $$\tau(\lambda)\otimes Y'= ((P_{\lambda}\circ\pi_{0,0})\otimes 
 Id)({\it\Delta}(Y)-\tau(\lambda)\otimes Y).$$  Then by (7.2),
 $Y+Y'=0$.  Hence, we obtain relations (iii)-(vi) by 
 evaluating  $$((P_{\lambda}\circ\pi_{0,0})\otimes 
 Id)({\it\Delta}(Y)-\tau(\lambda)\otimes Y)$$ to find $Y'$.

 By the Lemma 7.1 and [L4, (6.5) and (7.15)], we have $$\eqalign{&{\it\Delta}(B_i)=
 {\it\Delta}(y_it_i+d_i[(\adr Z_i)t_{p(i)}^{-1}x_{p(i)}]t_i+s_it_i)\cr&=y_it_i\otimes 
 1+t_i\otimes B_i+t_{p(i)}^{-1}x_{p(i)}t_i\otimes d_i[(\adr 
 Z_i)t_{p(i)}^{-2}]t_{p(i)}t_i\cr&+
 \sum_{\gamma>\alpha_i}G^+_{\gamma}t_i\otimes B.\cr}\leqno{(7.3)}$$  Note that the 
 $s_i$ are hidden in the above form of the coproduct.   This fact 
is the key reason   
  that the relations are independent of the choice of $s_i$. (See 
  also [L2, Lemma 5.7].)
 
The case  when  $\alpha_i$ is in
$\pi_{\Theta}$ follows directly from [L4, 7.20].  Indeed, checking the 
possibilities,     $a_{ij}$ must 
equal $0,-1,$ or 
$-2$ in this case.  Furthermore
$\delta_{p'(i),j}=\delta_{i,p'(j)}=\delta_{i,p'(i)}=0$ and thus (iii)-(v)
look just like 
 the quantum Serre relations between $y_it_i$ and $y_jt_j$.  This 
 is exactly what is explained in  the paragraph concerning [L4, 7.20]. Thus for the remainder of the proof,
 we assume  that $\alpha_i\notin\pi_{\Theta}$.
 
 Suppose that $\alpha_j\notin\pi_{\Theta}$.  By (i), it follows   $B_j$ commutes with 
each  $t_ix_i^{-1}$, $\alpha_i\in\pi_{\Theta}$.  Hence $B_j$ commutes 
with every element of ${\cal M}\cap G^+$.  In particular,  $[(\adr 
 Z_i)t_{p(i)}^{-2}] t_{p(i)}^2$, which is an element of $ {\cal 
 M}\cap 
 G^+$ for each $\alpha_i\notin\pi_{\Theta}$, commutes with $B_j$ for all 
 $\alpha_j\notin\pi_{\Theta}$.
 
 \medskip
 \noindent
 {\bf Case 1:} $a_{ij}=0$.    
 In this case, $Y=B_iB_j-B_jB_i$ and $\lambda=\alpha_i+\alpha_j$.
 Note that this expression is symmetric in $i$ and $j$ up to a 
 negative sign.  In particular, we may assume that $\alpha_j$ is also 
 not in $\pi_{\Theta}$.
 Expression (7.3) implies that  and $$\eqalign{&
 ((P_{\lambda}\circ\pi_{0,0})\otimes 
 Id)({\it\Delta}(B_iB_j)-\tau(\lambda)\otimes 
 B_iB_j)\cr &=(P_{\lambda}\circ\pi_{0,0})([t_{p(i)}^{-1}x_{p(i)}t_{i}y_{j}
 t_{j}
 \otimes d_i [(\adr Z_i)t_{p(i)}^{-2}]t_{p(i)}t_i).\cr}$$ 
 Since $t_{p(i)}^{-1}x_{p(i)}t_{i}y_{j}
 t_{j}=\delta_{p(i),j}(q_i-q_i^{-1})^{-1}t_i(t_j-t_j^{-1})$, it 
 follows that $$(P_{\lambda}\circ\pi_{0,0})(t_{p(i)}^{-1}x_{p(i)}t_{i}y_{j}
 t_{j})=\delta_{p(i),j}(q_i-q_i^{-1})^{-1}.$$  Thus
 $$\eqalign{&
 ((P_{\lambda}\circ\pi_{0,0})\otimes 
 Id)({\it\Delta}(B_iB_j)-\tau(\lambda)\otimes 
 B_iB_j)\cr &= \delta_{p(i),j}(q_i-q_i^{-1})^{-1}t_it_j\otimes 
 (d_i 
 [(\adr Z_i)t_{p(i)}^{-2}]t_{p(i)}
 t_i.\cr}$$  The same argument shows that 
 $$\eqalign{&
 ((P_{\lambda}\circ\pi_{0,0})\otimes 
 Id)({\it\Delta}(B_jB_i)-\tau(\lambda)\otimes 
 B_jB_i)\cr &= \delta_{p(j),i}(q_i-q_i^{-1})^{-1}t_it_j\otimes 
 (d_j 
 [(\adr Z_j)t_{p(j)}^{-2}]t_{p(j)}
 t_j.\cr}$$  Relation (iii) follows from the fact that 
 $\delta_{p(j),i}=\delta_{p(i),j}=\delta_{p'(i),j}$.

 \medskip
 \noindent
 {\bf Case 2:} $a_{ij}=-1$.  In this case $$Y=
 B_i^2B_j-(q_i+q_i^{-1})B_iB_jB_i+B_jB_i^2.$$  Recall that if 
 $p(i)=i$, then $d_i=1$. Furthermore, checking 
 the possibilities, we have $i=p(j)$ if and only if 
 $\Theta(\alpha_i)=-\alpha_{p(i)}$ and 
 $\Theta(\alpha_j)=-\alpha_{p(j)}$.   Thus 
 $$
 (P_{\lambda}\circ\pi_{0,0})\otimes 
 Id)({\it\Delta}(Y)-\tau(\lambda)\otimes 
 Y) = 
 (P_{\lambda}\circ\pi_{0,0})
 (\delta_{i,p(i)}(G_1)+\delta_{i,p(j)}(G_2))$$
 where $$\eqalign{&G_1=t_i^{-1}x_it_iy_it_it_j\otimes  
 [(\adr Z_i)t_i^{-2}]t_i
 t_iB_j\cr&-(q_i+q_i^{-1})t_i^{-1}x_it_it_jy_it_i\otimes  
 [(\adr Z_i)t_i^{-2}]t_i
 t_iB_j\cr&+t_jt_i^{-1}x_it_iy_it_i\otimes B_j[(\adr 
 Z_i)t_i^{-2}]t_it_i\cr}$$
 and
 $$\eqalign{&G_2=(t_{j}^{-1}x_{j}t_it_iy_jt_j\otimes 
 d_it_{j}^{-1}t_iB_i+t_it_{j}^{-1}x_{j}t_iy_jt_j\otimes
B_i d_it_{j}^{-1}t_i\cr& -(q_i+q_i^{-1})(t_{j}^{-1}x_{j}t_iy_jt_jt_i\otimes 
 d_it_{j}^{-1}t_iB_i+t_it_i^{-1}x_it_jy_it_i\otimes d_jB_it_i^{-1}t_j)\cr
&+(t_{i}^{-1}x_{i}t_jt_iy_it_i\otimes 
 d_jt_{i}^{-1}t_jB_i+t_{i}^{-1}x_{i}t_jy_it_it_i\otimes
 d_jt_{i}^{-1}t_jB_i).\cr }$$
Suppose that $\beta_1,\beta_2,\beta_3$ are elements of $Q(\pi)$ such 
that $\beta_1+\beta_2+\beta_3=2\alpha_i$. Recall that 
$\lambda=2\alpha_i+\alpha_j$. Note that $$(P_{\lambda}\circ 
\pi_{0,0})(\tau(\beta_1)x_j
  \tau(\beta_1)y_j\tau(\beta_3)
  =(q_i-q_i^{-1})^{-1}q^{-(\beta_1,\alpha_j)}t_i^2t_j.$$
 Hence $$\eqalign{ (P_{\lambda}\circ\pi_{0,0})(G_1)&=
 (q_i-q_i^{-1})[t_i^2t_j(q_i^{-2}-(q_i+q_i^{-1})q_i^{-1}+q_i^{-2})\otimes 
 [(\adr Z_i)t_i^{-2}]t_i^2B_j
 \cr
 &=
 -q_i^{-1}t_i^2t_j\otimes [(\adr Z_i)t_i^{-2}]t_i^2B_j.\cr}$$
 A similar computation  shows that
 $$\eqalign{ &(P_{\lambda}\circ\pi_{0,0})(G_2)=
 (q_i-q_i^{-1})[t_i^2t_j(q_i^{2}q_i^{-3}+q_i-(q_i+q_i^{-1})q_iq_i^{-3})\otimes 
 d_iB_it_j^{-1}t_i\cr
 &+t_i^2t_j(-(q_i+q_i^{-1})q_i+q_i^{-1}q_i^3+q_i  q_i^{3})\otimes
 d_jB_it_i^{-1}t_j]\cr
 &=t_i^2t_j\otimes 
 [q_i^{-1}(q_i+q_i^{-1})d_iB_it_j^{-1}t_i
 +q_i^{2}(q_i+q_i^{-1})d_jB_it_i^{-1}t_j].\cr}$$
 
 \medskip
 \noindent
 {\bf Case 3}: $a_{ij}=-2$.  Checking the possibilities, we see that 
 $\Theta(\alpha_i)=-\alpha_i$ and $\Theta(\alpha_j)=-\alpha_j$.  
 Furthermore, ${\cal S}$ is equal to $\{\alpha_j\}$.  
 Thus $B_i=y_it_i+q_i^{-2}x_i$ and $B_j=y_jt_j+q_j^{-2}x_j+s_jt_j$.  
 It follows that ${\it\Delta}(B_i)=B_i\otimes 
 1+t_i\otimes B_i$ and ${\it\Delta}(B_j)= (y_jt_j+q_j^{-2}x_j)\otimes 
 1+t_j\otimes B_j$.  So 
$$
\eqalign{ &(P_{\lambda}\circ\pi_{0,0})\otimes 
 Id)({\it\Delta}(Y)-\tau(\lambda)\otimes 
 Y)\cr&=(P_{\lambda}\circ\pi_{0,0})\otimes 
 Id)(H_1\otimes B_iB_j+H_2\otimes B_jB_i)\cr}$$
 where $$\eqalign{H_1&=B_i^2t_it_j+B_it_iB_it_j+t_iB_iB_it_j
 \cr&+(q_i^2+1+q_i^{-2})[-t_iB_it_jB_i-B_it_it_jB_i
 +t_it_jB_iB_i]
 \cr}$$ and 
$$\eqalign{H_2&= (q_i^2+1+q_i^{-2})(-B_iB_it_jt_i+
  B_it_jB_it_i+B_it_jt_iB_i)
 \cr&-t_jB_iB_it_i-t_jB_it_iB_i-t_jt_iB_iB_i.\cr}$$

 Let $\tau_1=\tau(\beta_1),\tau_2=\tau(\beta_2)$ and $\tau_3=\tau(\beta_3)$
 where $\beta_1,\beta_2, $ and $\beta_3$ are three  elements in $Q(\pi)$ such that 
 $\beta_1+\beta_2+\beta_3=\alpha_i+\alpha_j$.
  Relation (iv) now follows from the fact that   
 $$(P_{\lambda}\circ\pi_{0,0})\tau_1
 B_i\tau_2B_i\tau_3=(P_{\lambda}\circ\pi_{0,0})
(\tau_1q_i^{-2}x_i\tau_2y_it_i\tau_3)
 =q^{-(\beta_2,\alpha_i)}(q_i-q_i^{-1})^{-1}t_i^3t_j.$$
 
\medskip
 \noindent
 {\bf Case 4}:  $a_{ij}=-3$.   Again, checking the possibilities, we 
 must have that $\Theta(\alpha_i)=-\alpha_i$,
 $\Theta(\alpha_j)=-\alpha_j$, and ${\cal S}$ is empty. 
 Thus $B_k=y_kt_k+q_k^{-2}x_k$ and ${\it\Delta}(B_k)=B_k\otimes 
 1+t_k\otimes B_k$ for $k=1,2$. The argument follows in a similar, 
 although lengthier, manner to Case 3.
 $\Box$ 

\bigskip
\centerline{REFERENCES}
\medskip

{\medskip
\noindent[A] S. Araki, On root systems and an infinitesimal 
classification of irreducible symmetric spaces, {\it Journal of 
Mathematics, Osaka City University} {\bf 13} (1962), no. 1, 1-34.

\medskip
\noindent [Di] M.S. Dijkhuizen, Some remarks on the construction
of quantum symmetric spaces, In: {\it Representations of Lie
Groups, Lie Algebras and Their Quantum Analogues, } Acta Appl.
Math. {\bf 44} (1996), no. 1-2, 59-80.

\medskip
\noindent [DN] M.S. Dijkhuizen and M. Noumi, A family of quantum
projective spaces and related $q$-hypergeometric orthogonal
polynomials, {\it Trans. Amer. Math. Soc.} {\bf 350} (1998), no.
8, 3269-3296.

\medskip
\noindent [D] J. Dixmier, {\it Alg\`ebres Enveloppantes}, Cahiers
Scientifiques, XXXVII,

\noindent Gauthier-Villars, Paris (1974).

\medskip
\noindent
[GK] A.M. Gavrilik and A. U. Klimyk, 
$q$-deformed orthogonal and pseudo-orthogonal algebras and their representations,
{\it Lett. Math. Phys.}{\bf 21} (1991), no. 3, 215-220.

\medskip
\noindent
[GI] A.M. Gavrilik and N.Z. Iorgov, $q$-Deformed algebras $U_q({\bf 
so}\ n)$ and their representations, q-alg/9709036.

\medskip
\noindent
[GIK] A.M. Gavrilik and N.Z. Iorgov,
Representations of the nonstandard algebras 
$U\sb q({\rm so}(n))$ and $U\sb q({\rm so}(n-1,1))$ in Gelfand-Tsetlin basis,
{\it International Symposium on Mathematical and Theoretical Physics} (Kyiv, 1997),
Ukra•n. Fi z. Zh.{\bf 43} (1998), no. 6-7, 791-797. 

\medskip
\noindent
[HKP] 
M. Havl'\v cek, A.U. Klimyk, and S. Po\v sta, 
Representations of the $q$-deformed algebra $U\sb q\prime({\rm so}\sb 4)$. 
J. Math. Phys.{\bf 42} (2001), no. 11, 5389-5416. 

\medskip
\noindent
[He] S. Helgason, {\it Differential Geometry, Lie Groups, and Symmetric 
Spaces,}
 Pure and Applied Mathematics {\bf 80}, Academic Press,  New York  (1978).
 
 \medskip
 \noindent
 [He2] S. Helgason, {\it Geometric Analysis on Symmetric Spaces,}  
 Mathematical surveys and monographs, American Mathematical Society,  
{\bf 39} (1994).
 
 \medskip
\noindent [H] J.E. Humphreys, {\it Introduction to Lie Algebras
and Representation Theory}, Springer-Verlag, New York
(1972).

\medskip
\noindent [J] N. Jacobson, {\it Basic Algebra II,}
W. H. Freeman and Company, San Francisco  (1980).

\medskip
\noindent [Jo] A. Joseph, {\it Quantum Groups and Their Primitive
Ideals}, Springer-Verlag, New York (1995).

\medskip
\noindent [K1] T.H.  Koornwinder, Askey-Wilson polynomials as zonal 
spherical functions on the SU(2) quantum group, {\it SIAM Journal on 
Mathematical Analysis}, {\bf 24} (1993), no. 3, 795-813.

\medskip
\noindent [K2] T.H. Koornwinder, Askey-Wilson polynomials for root 
systems of type BC, In: {\it Hypergeometric functions on domains of 
positivity, Jack polynomials and applications}, Contemporary 
Mathematics {\bf 138} (1992) 189-204. 

\medskip \noindent [KS] A.  Klimyk and K. Schmudgen,  {\it Quantum
Groups and Their Representations},  Texts and Monographs in
Physics, Springer-Verlag, Berlin (1997).

\medskip
\noindent
[Kn] A. W.  Knapp, {\it Lie Groups Beyond an Introduction,}
 Progress in Math. {\bf 140}, Birkh\"auser,  Boston
 (1996).

\medskip
\noindent [L1] G. Letzter, Subalgebras which appear in quantum
Iwasawa decompositions, {\it Canadian Journal of Math.}
  {\bf 49} (1997), no. 6, 1206-1223.
   
\medskip \noindent [L2] G. Letzter,  Symmetric pairs for quantized
enveloping algebras, {\it Journal of Algebra} {\bf 220} (1999), no. 2,
729-767.

\medskip \noindent [L3] G. Letzter, Harish-Chandra modules for
quantum symmetric pairs, {\it Representation Theory, An Electronic
Journal of the AMS} {\bf 4} (1999)
64-96.

\medskip \noindent [L4] G. Letzter, Coideal subalgebras and quantum 
symmetric pairs, MSRI volume, 1999 Hopf Algebra Workshop, to appear,
(math.QA/0103228).

\medskip\noindent
[Ma] I.G. Macdonald, Orthogonal polynomials associated with root 
systems, {\it S\'eminaire Lotharingien de Combinatoire}, {\bf 45} 
(2001), Article B45a.
 
 \medskip\noindent
[M] S. Montgomery,  {\it Hopf Algebras and Their Actions on Rings,}
 CBMS Regional Conference Series in Mathematics {\bf 82},
  American Mathematical Society, Providence (1993).

\medskip
\noindent [NS] M. Noumi and T. Sugitani, Quantum symmetric spaces
and related q-orthogonal polynomials, In: {\it Group Theoretical
Methods in Physics (ICGTMP)} (Toyonaka, Japan, 1994), World
Sci.~Publishing, River Edge, N.J. (1995) 28-40.

\medskip
\noindent
[OV] A. L. Onishchik and E. B. Vinberg, {\it Lie Groups and Lie 
Algebras  III:
Structure of Lie Groups and Lie Algebras,}  Springer-Verlag, Berlin (1994).}

\end{document}